\input amstex
\documentstyle{amsppt}

\topmatter
\title
        Hochschild (co)homology of Hopf crossed products
\endtitle

\author
       Jorge A. Guccione and Juan J. Guccione
\endauthor

\address
     Jorge Alberto Guccione, Departamento de Matem\'atica, Facultad de
     Ciencias Exactas y Naturales, Pabell\'on 1 - Ciudad Universitaria,
     (1428) Buenos Aires, Argentina.
\endaddress

\email
     vander\@dm.uba.ar
\endemail

\address
     Juan Jos\'e Guccione, Departamento de Matem\'atica, Facultad de
     Ciencias Exactas y Naturales, Pabell\'on 1 - Ciudad Universitaria,
     (1428) Buenos Aires, Argentina.
\endaddress

\email
     jjgucci\@dm.uba.ar
\endemail

\abstract
For a general crossed product $E=A\#_f H$, of an algebra $A$ by a Hopf
algebra $H$, we obtain complexes simpler than the canonical ones, giving
the Hochschild homology and cohomology of $E$. These complexes are equipped
with natural filtrations. The spectral sequences associated to them is a
natural generalization of the one obtained in \cite{H-S} by the direct
method. We also get that if the $2$-cocycle $f$ takes its values in a
separable subalgebra of $A$, then the Hochschild (co)homology of $E$ with
coefficients in $M$ is the (co)homology of $H$ with coefficients in a
(co)chain complex.

\endabstract


\thanks
Supported by UBACYT 01/TW79 and CONICET
\endthanks


\endtopmatter

\document

\def \circ{\,}
\def \ni{\noindent}
\def \ot{\otimes}
\def \ov{\overline}
\def \wt{\widetilde}
\def \wh{\widehat}
\def \sup{\supseteq}
\def \sub{\subseteq}
\def \ba{\bold a}
\def \bb{\bold b}

\def \bbf{\bold f}
\def \bg{\bold g}
\def \bh{\bold h}

\def \bx{\bold x}

\def \bQ{\Bbb Q}
\def \bby{\bold y}

\def \fh{\frak h}

\def \fg{\frak g}

\def \de{\delta}
\def \ep{\epsilon}

\def \si{\sigma}
\def \om{\omega}
\def \Om{\Omega}

\def \B{\operatorname{B}}
\def \C{\operatorname{C}}
\def \Ext{\operatorname{Ext}}
\def \gl.dim{\operatorname{gl.dim}}
\def \h{\operatorname{h}}
\def \H{\operatorname{H}}
\def \Hom{\operatorname{Hom}}

\def \op{\operatorname{op}}

\def \Tor{\operatorname{Tor}}

\head  Introduction \endhead
Let $G$ be a group, $S = \bigoplus S_g$ a strongly $G$-graded algebra and
$V$ an $S$-bimodule. In \cite{L} was shown that there is a convergent
spectral sequence
$$
E^2_{rs} = \H_r(G,\H_s(S_e,V)) \Rightarrow \H_{r+s}(S,V),
$$
where $e$ denotes the identity of $G$. In \cite{S} was shown that this
result remains valid for $H$-Galois extensions (in his paper the author
deals with both the homology and the cohomology of these algebras). An
important particular type of $H$-Galois extensions are the crossed products
with convolution invertible cocycle $E=A\#_f H$, of an algebra $A$ by a
Hopf algebra $H$ (for the definition see Section one). The purpose of our
paper is to construct complexes simpler than the canonical ones, given the
Hochschild (co)homology of $E$ with coefficients in an arbitrary
$E$-bimodule. These complexes are equipped with canonical filtrations. We
show that the spectral sequences associated to them coincide with the ones
obtained using a natural generalization of the direct method introduced in
\cite{H-S}, and with the ones constructed in \cite{S} (when these are
specialize to crossed products). In the case of group extensions these
results were proved in \cite{E} and \cite{B}.

\smallskip

This paper is organized as follows: in Section~1 a resolution $(X_*,d_*)$
of a crossed product $E=A\#_f H$ is given. To accomplish this construction
we do not use the fact that the cocycle is convolution invertible.
Moreover, we give a recursive construction of morphisms $\phi_*\: (X_*,d_*)
\to (E\ot \ov{E}^*\ot E, b'_*)$ and $\psi_*\: (E\ot \ov{E}^*\ot E, b'_*)
\to (X_*,d_*)$, where $(E\ot \ov{E}^*\ot E, b'_*)$ is the normalized
Hochschild resolution, such that $\psi_*\circ\phi_* = id$ and we show that
$\phi_*\circ\psi_*$ is homotopically equivalent to the identity map.
Consequently our resolution is a direct sum of the normalized Hochschild
resolution. We also recursively construct an homotopy $\phi_*\circ \psi_*
@>\om_{*+1}>> id_*$. Both, the canonical normalized resolution and
$(X_*,d_*)$ are equipped with natural filtrations, which are preserved by
the maps $\phi_*$, $\psi_*$ and $\om_{*+1}$.

In Section~2, for an $E$-bimodule $M$, we get complexes $\wh{X}_*(E,M)$
and $\wh{X}^*(E,M)$, giving the Hochschild homology and cohomology of $E$
with coefficients in $M$ respectively. The filtration of $(X_*,d_*)$
induces filtrations on $\wh{X}_*(E,M)$ and $\wh{X}^*(E,M)$. So, we
obtain converging spectral sequences $E^1_{rs} = \H_r(A,M\ot \ov{H}^s)
\Rightarrow \H_{r+s}(E,M)$ and $E_1^{rs} = \H^r(A,\Hom_k(\ov{H}^s,M))
\Rightarrow \H^{r+s}(E,M)$. Using the results of Section~1, we get that
these spectral sequences are the ones associated to suitable filtrations of
the Hochschild normalized chain and cochain complexes $(M\ot \ov{E}^*,b_*)$
and $(\Hom_k(\ov{E}^*,M),b^*)$. This allows us to give very simple proofs of
the main results of \cite{H-S} and \cite{G}.

In Section~3, we show that, if the cocycle is convolution invertible, then
the complexes $\wh{X}_*(E,M)$ and $\wh{X}^*(E,M)$ are isomorphic to simpler
complexes $\ov{X}_*(E,M)$ and $\ov{X}^*(E,M)$ respectively. Then, we
compute the term $E^2_{rs}$ and $E_2^{rs}$ of the spectral sequences
obtained in Section~2. Moreover, using the above mentioned filtrations, we
prove that if the $2$-cocycle $f$ takes its values in a separable
subalgebra of $A$, then the Hochschild (co)homology of $E$ with
coefficients in $M$ is the (co)homology of $H$ with coefficients in a
(co)chain complex. Finally, as an application we obtain some results about
the $\Tor_*^E$ and $\Ext^*_E$ functors and an upper bound for the global
dimension of $E$ (for group crossed products this bound was obtained in
\cite{A-R}).

In addition to the direct method developed in \cite{H-S}, there are another
two classical methods to obtain spectral sequences converging to
$\H_*(E,M)$ and with $E^2$-term $\H_*(H,\H_*(A,M))$. Namely the
Cartan-Leray and the Grothendieck spectral sequences of a crossed product.
In Section~4, we recall these constructions and we prove that these
spectral sequences are isomorphic to the one obtained in Section~2. This
generalizes the main results of \cite{B}.

In a first appendix we give a method to construct (under suitable
hypothesis) a projective resolution of the $k$-algebra $E$ as $E^e = E\ot
E^{\op}$-bimodule, simpler than the canonical one of Hochschild. This
method, which can be considered as a variant of the perturbation lemma, is
used to prove the main result of Section~1. The boundary maps of the
resolution $(X_*,d_*)$ are recursively defined in Section~1. In a second
appendix we give closed formulas for these maps.

\head 1. A resolution for a crossed product \endhead
Let $A$ be a $k$-algebra and $H$ a Hopf algebra. We will use the Sweedler
notation $\Delta(h) = h^{(1)} \ot h^{(2)}$, with the summation understood
and superindices instead of subindices. Recall some definitions of
\cite{B-C-M} and \cite{D-T}. A {\it weak action} of $H$ on $A$ is a
bilinear map $(h,a)\mapsto a^h$ from $H\times A$ to $A$ such that, for
$h\in H$, $a,b\in A$

\medskip

\itemitem{1)} $(ab)^h = a^{h^{(1)}}b^{h^{(2)}}$,

\smallskip

\itemitem{2)} $1^h = \ep(h)1$,

\smallskip

\itemitem{3)} $a^1 = a$.

\medskip

Let $A$ be a $k$-algebra and $H$ a Hopf algebra with a weak action on $A$.
Given a $k$-linear map $f\:H\ot H\to A$, let $A\#_f H$ be the $k$-algebra (in
general non associative and without $1$) with underlying vector space $A\ot
H$ and multiplication map
$$
(a\ot h)(b \ot l) =  a b^{h^{(1)}}f(h^{(2)},l^{(1)}) \ot h^{(3)}l^{(2)},
$$
for all $a,b\in A$, $h,l \in H$. The element $a\ot h$ of $A\#_f H$ will
usually be written $a\# h$ to remind us $H$ is weakly acting on $A$. The
algebra $A\#_f H$ is called a {\it crossed product} if it is associative with
$1\# 1$ as identity element. It is easy to check that this happens if and
only if $f$ and the weak action satisfy the following conditions:

\medskip

\itemitem{i)} (Normality of $f$) for all $h\in H$, we have $f(h,1) = f(1,h)
= \ep(h)1_A$,

\smallskip

\itemitem{ii)} (Cocycle condition) for all $h,l,m\in H$, we have
$$
f\bigl(l^{(1)},m^{(1)}\bigr)^{h^{(1)}} f\bigl(h^{(2)}, l^{(2)}m^{(2)}\bigr)
=  f\bigl(h^{(1)},l^{(1)}\bigr) f\bigl(h^{(2)}l^{(2)},m\bigr),
$$

\smallskip

\itemitem{iii)} (Twisted module condition) for all $h,l\in H$, $a\in A$
we have
$$
\bigl(a^{l^{(1)}}\bigr)^{h^{(1)}} f\bigl(h^{(2)},l^{(2)} \bigr) =
f\bigl(h^{(1)},l^{(1)}\bigr)a^{h^{(2)}l^{(2)}}.
$$

\medskip

In this section we obtain a resolution $(X_*,d_*)$ of a crossed product $E
= A\#_f H$ as an $E$-bimodule, which is simpler than the canonical one of
Hochschild. To begin, we fix some notations:

\medskip

\item{1)} For each $k$-algebra $B$, we put $\ov{B} = B/k$. Moreover, given
$b\in B$ we also let $b$ denote the class of $b$ in $\ov{B}$. 

\smallskip

\item{2)} We write $B^l = B\ot\cdots\ot B$, $\ov{B}^l = \ov{B}\ot\cdots\ot
\ov{B}$ ($l$ times) and $\B_l (B)= B\ot \ov{B}^l\ot B$, for each
natural number $l$.

\smallskip

\item{3)} Given $a_0\ot\cdots\ot a_r \in A^{r+1}$ and $0\le i<j\le r$, we
write $\ba_{ij} = a_i\ot\cdots\ot a_j \in A^{j-i+1}$.

\smallskip

\item{4)} Given $h_0\ot\cdots\ot h_s \in H^{s+1}$ and $0\le i<j\le s$, we
write $\bh_{ij} = h_i\ot\cdots\ot h_j$ and $\fh_{ij} = h_i\cdots h_j \in
H$.

\smallskip

\item{5)} Given $\bh = h_0\ot\cdots\ot h_s\in H^{s+1}$, we let
$\bh^{(1)}\ot \bh^{(2)}$ denote the comultiplication of $\bh$ in $H^{s+1}$.
So, $\bh^{(1)}\ot\bh^{(2)} = (h_0^{(1)}\ot \cdots\ot h_s^{(1)}) \ot
(h_0^{(2)}\ot\cdots\ot h_s^{(2)})$.

\smallskip

\item{6)} Given $a\in A$, $\ba = a_1\ot\cdots\ot a_r\in A^r$ and $\bh =
h_0\ot\cdots\ot h_s\in H^{s+1}$, we write $a^{\ov{\bh}} =
(\dots(((a^{h_s})^{h_{s-1}})^{ h_{s-2}})^{ h_{s-3}}\dots)^{h_0}$ and
$\ba^{\ov{\bh}} = a_1^{\ov{\bh_{0s}^{(1)}}}\ot\cdots\ot
a_r^{\ov{\bh_{0s}^{(r)}}}$.

\specialhead 1.1. The resolution $(X_*,d_*)$\endspecialhead
Let $Y_s = E\ot \ov{H}^s \ot H$ ($s\ge 0$) and
$X_{rs} = E\ot \ov{H}^s \ot \ov{A}^r \ot E$ ($r,s \ge 0$). The groups
$X_{rs}$ are $E$-bimodules in an obvious way and the groups $Y_s$ are
$E$-bimodules via the left canonical action and the right action
$$
(a_0\ot\bh)(a\# h) =  a_0a^{\ov{ \bh^{(1)}}}
f(h_{s+1}^{(2)},h^{(1)})^{\ov{\bh_{0s}^{(2)}}} \ot \bigl(\bh_{0s}^{(3)}
\ot h_{s+1}^{(3)} h^{(2)}\bigr),
$$
where $\bh = h_0\ot \cdots\ot h_{s+1}$. Let us consider the diagram of
$E$-bimodules and $E$-bimodule maps
$$
\CD
\vdots\\
@VV\partial_2V \\
Y_1 @<\mu_1<<  X_{01} @<d^0_{11}<<  X_{11} @<d^0_{21}<< \dots \\
@VV\partial_1V \\
Y_0 @<\mu_0<<  X_{00} @<d^0_{10}<<  X_{10} @<d^0_{20}<< \dots,
\endCD
$$
where $\mu_*\: X_{0*} \to Y_*$, $d^0_{**}\: X_{**} \to X_{*-1,*}$ and
$\partial_*\: Y_* \to Y_{*-1}$ are defined by:
$$
\align
&\mu_s(a_0\ot\bh_{0s}\ot a_1\ot h_{s+1}) = a_0a_1^{\ov{\bh_{0s}^{(1)}}}\ot
\bh_{0s}^{(2)} \ot h_{s+1},\\
&d^0_{rs}(a_0\ot\bh_{0s}\ot\ba\ot h_{s+1}) = a_0a_1^{\ov{ \bh_{0s}
^{(1)}}}\ot\bh_{0s}^{(2)}\ot \ba_{2,r+1}\ot h_{s+1}\\
&\phantom{d^0_{rs}(a_0\ot\bh_{0s}\ot \ba\ot h_{s+1})}
+\sum_{i=1}^r (-1)^i a_0\ot\bh_{0s}\ot \ba_{1,i-1}\ot a_ia_{i+1}\ot
\ba_{i+2,r+1}\ot h_{s+1},\\
&\partial_s(a\ot\bh) = \sum_{i=0}^s (-1)^{i+1} a f(h_i^{(1)},
h_{i+1}^{(1)})^{\ov{\bh^{(1)}_{i-1}}} \ot \bh_{0,i-1}^{(2)} \ot h_i^{(2)}
h_{i+1}^{(2)} \ot \bh_{i+2,s+1},
\endalign
$$
where $\ba = a_1\ot\cdots\ot a_{r+1}$ and $\bh = h_0\ot \cdots \ot
h_{s+1}$. We have left $E$-module maps $\si^0_{0*} \:Y_* \to X_{0*}$ and
$\si^0_{*+1,*} \:X_{**} \to X_{*+1,*}$, given by
$\si^0_{r+1,s}(a_0\ot\bh_{0s}\ot \ba\ot h_{s+1})= (-1)^{r+1} a_0\ot
\bh_{0s}\ot \ba\ot 1\# h_{s+1}$ for $r\ge -1$. Clearly $(Y_*,\partial_*)$
is a complex and $\si^0_{*+1,s}$ is a contracting homotopy of
$$
Y_s @<\mu_s <<  X_{0s} @<d^0_{1s}<< X_{1s} @<d^0_{2s}<< X_{2s} @<d^0_{2s}<<
X_{3s} @<d^0_{3s}<< X_{4s} @<d^0_{4s}<< X_{5s} @<d^0_{5s}<<\dots.
$$
So, we are in the situation considered in Appendix~A. We define
$E$-bimodule maps $d^l_{rs}\: X_{rs} \to X_{r+l-1,s-l}$ ($r\ge 0$ and $1\le
l\le s$) recursively, by:
$$
d^l_{rs}(\bx) = \cases
- \si_{0,s-1}^0\circ\partial_s\circ\mu_s(\bx) &\text{if $r=0$ and $l=1$,}\\
- \sum_{j=1}^{l-1} \si_{l-1,s-l}^0\circ d^{l-j}_{j-1,s-j}\circ d^j_{0s}
(\bx) & \text{if $r=0$ and $1<l\le s$,}\\
- \sum_{j=0}^{l-1} \si_{r+l-1,s-l}^0\circ d^{l-j}_{r+j-1,s-j}\circ
d^j_{rs}(\bx) &\text{if $r>0$,}
\endcases
$$
for $\bx \in k\ot\ov{H}^s\ot \ov{A}^r\ot k$.

\proclaim{Theorem 1.1.1} There is a relative projective resolution
$$
E @<-\mu<< X_0  @<d_1<< X_1  @<d_2<< X_2  @<d_3<< X_3 @<d_4<< X_4 @<d_5<<
X_5 @<d_6<< X_6@<d_7<<\dots, \tag 1
$$
where $\displaystyle{X_n = \bigoplus_{r+s=n} X_{rs}}$, $\mu$ is the
multiplication map  and $\displaystyle{d_n = \sum_{r+s=n\atop r+l> 0}
\sum^s_{l=0} d^l_{rs}}$.
\endproclaim

\demo{Proof} Let $\wt{\mu}\:Y_0\to E$ be the map $\wt{\mu}(a\ot(h_0\ot
h_1)) = -af(h_0^{(1)},h_1^{(1)})\# h_0^{(2)}h_1^{(2)}$. The complex of
$E$-bimodules
$$
E @<\wt{\mu}<< Y_0  @<\partial_1<< Y_1  @<\partial_2 << Y_2 @<\partial_3<<
Y_3 @<\partial_4<< Y_4  @<\partial_5 << Y_5 @<\partial_6<< Y_6
@<\partial_7<<\dots
$$
is contractible as a complex of left $E$-modules. A chain contracting
homotopy $\si_0^{-1}\: E \to Y_0$ and $\si^{-1}_{s+1}\: Y_s\to Y_{s+1}$
($s\ge 0$) is given by $\si^{-1}_{s+1}(x) = (-1)^s x\ot 1_H$. Hence, the
theorem follows from Corollary~A.2 of Appendix~A\qed

\enddemo

\remark{Remark 1.1.2}  Let $\si^l_{l,s-l}\: Y_s \to X_{l,s-l}$ and
$\si^l_{r+l+1,s-l}\: X_{rs} \to X_{r+l+1,s-l}$  be the maps recursively
defined by
$$
\si^l_{r+l+1,s-l} = - \sum_{i=0}^{l-1} \si^0_{r+l+1,s-l} \circ d^{l-i}
_{r+i+1,s-i} \circ\si^i_{r+i+1,s-i}\quad\text{($0<l\le s$ and $r\ge -1$)}.
$$
We will prove, in Corollary~A.2, that the family $\ov{\si}_0\: E \to X_0$,
$\ov{\si}_{n+1}\:X_n \to X_{n+1}$, defined by $\ov{\si}_0 =
\si_{00}^0 \circ \si_0^{-1}$ and
$$
\ov{\si}_{n+1} = - \sum_{l=0}^{n+1} \si_{l,n-l+1}^l \circ \si_{n+1}^{-1}
\circ \mu_n + \sum_{r+s=n} \sum_{l=0}^s \si_{r+l+1,s-l}^l\quad \text{($n\ge
0$)},
$$
is a contracting homotopy of the resolution $(1)$ introduced in
Theorem~1.1.1.
\endremark

\proclaim{Theorem 1.1.3} Let $\bx= a_0\ot \bh\ot \ba\ot 1_E$, with $\ba =
a_1\ot\cdots\ot a_r\in \ov{A}^r$ and $\bh = h_0\ot\cdots\ot h_s\in H\ot
\ov{H}^s$. We have:

\smallskip

\item{1)} $d^1_{rs}$ is the map given by
$$
\align
d^1_{rs}(\bx) &= 
\sum_{i=0}^{s-1} (-1)^{i+r} a_0 f(h_i^{(1)},h_{i+1}^{(1)})^
{\ov{\bh_{0,i-1}^{(1)}}}\ot \bh_{0,i-1}^{(2)}\ot \h_i^{(2)} h_{i+1}^{(2)} \ot
\bh_{i+2,s} \ot \ba\ot 1_E\\
& + (-1)^{r+s} a_0\ot \bh_{0,s-1} \ot \ba^{h_s^{(1)}}\ot
1\# h_s^{(2)},
\endalign
$$

\smallskip

\item{2)} For each $l\ge 2$, there are maps $F_0^{(l)}\:\ov{H}^l \to
A^{l-1}$ and $F_r^{(l)}\:\ov{H}^l\ot \ov{A}^r \to A^{r+l-1}$ ($r\ge 1$),
whose image is included in the $k$-submodule of $A^{r+l-1}$ generated by all
the elementary tensors $a_1\ot\cdots\ot a_{r+l-1}$ with $l-1$ coordinates
in the image of $f$, such that for $2\le l\le s$,
$$
d^l_{rs}(\bx) = (-1)^{l(r+s)} a_0\ot \bh_{0,s-l}\ot F_r^{(l)}
(\bh_{s-l+1,s}^{(1)}\ot \ba) \ot 1\# \fh_{s-l+1,s}^{(2)},
$$
where $F_r^{(l)}(\bh_{s-l+1,s}^{(1)}\ot \ba) = F_0^{(l)}(\bh_{s-l+1,s})$
if $r=0$.
\endproclaim

\demo{Proof} The computation of $d_{rs}^1$ can be obtained easily by
induction on $r$, using that  $d_{0s}^1 = -\si_{0,s-1}^0\circ \partial_s
\circ \mu_s^0$ and $d_{rs}^1 = -\si_{r,s-1}^0\circ d_{r-1,s}^1 \circ
d_{rs}^0$ for $r\ge 1$. The assertion for $d_{rs}^l$, with $l\ge 2$,
follows easily by induction on $l$ and $r$, using the recursive definition
of $d_{rs}^l$\qed
\enddemo

In Appendix~B we will give more precise formulas for the maps $F_r^{(l)}$
completing the computation of the $d_{rs}^l$'s.

\specialhead 1.2. Comparison with the canonical resolution \endspecialhead
Let $(\B_* (E), b'_*)$ be the normalized Hochschild resolution of $E$. As
it is well known, the complex
$$
E @<\mu <<  E\ot E @<b'_1<< \B_1 (E) @<b'_2<< \B_2 (E) @<b'_3<<
\B_3 (E) @<b'_4<<\dots
$$
is contractible as a complex of left $E$-modules, with contracting homotopy
$\xi_n(\bx)\! = (-1)^n \bx\ot 1$. Let $\ov{\si}_*$ be the contracting
homotopy of $(1)$ introduced in Remark 1.1.2. Let $\phi_*\: (X_*,d_*) \to
(\B_* (E), b'_*)$ and $\psi_*\: (\B_*(E), b'_*) \to (X_*,d_*)$ be the
morphisms of $E$-bimodule complexes, recursively defined by $\phi_0 = id$,
$\psi_0 = id$, $\phi_{n+1}(\bx\ot 1) = \xi_{n+1} \circ \phi_n \circ
d_{n+1}(\bx\ot 1)$ and $\psi_{n+1}(\bby\ot 1) = \ov{\si}_{n+1} \circ \psi_n
\circ b'_{n+1}(\bby\ot 1)$.

\proclaim{Proposition 1.2.1} $\psi_*\circ\phi_* = id_*$ and $\phi_*\circ
\psi_*$ is homotopically equivalent to the identity map. An homotopy
$\phi_*\circ \psi_*  @>\om_{*+1}>> id_*$ is recursively defined by $\om_1 =
0$ and $\om_{n+1}(\bx) = \xi_{n+1}\circ (\phi_n \circ \psi_n - id -
\om_n \circ b'_n)(\bx)$, for $\bx\in E\ot\ov{E}^n\ot k$.

\endproclaim

\demo{Proof} We prove both assertions by induction. Let $U_n = \phi_n \circ
\psi_n - id_n$ and $T_n = U_n - \om_n\circ b'_n$. Assuming that $b'_n\circ
\om_n + \om_{n-1}\circ b'_{n-1} = U_{n-1}$, we get that on $E\ot \ov{E}^n
\ot k$,
$$
\align
b'_{n+1}\circ \om_{n+1} + \om_n\circ b'_n & = b'_{n+1}\circ \xi_{n+1}\circ
T_n+ \om_n\circ b'_n\\
&= T_n- \xi_n \circ b'_n\circ T_n + \om_n\circ b'_n\\
& = U_n - \xi_n \circ U_{n-1} \circ b'_n + \xi_n \circ b'_n\circ \om_n\circ
b'_n \\
&= U_n - \xi_n \circ U_{n-1}\circ b'_n + \xi_n \circ T_{n-1} \circ b'_n =
U_n.
\endalign
$$
Hence, $b'_{n+1}\circ \om_{n+1} + \om_n\circ b'_n = U_n$ on $\B_n(E)$. Next,
we prove that $\psi_*\circ \phi_* = id_*$. It is clear that $\psi_0
\circ\phi_0 = id_0$. Assume that $\psi_n\circ\phi_n = id_n$. Since
$\phi_{n+1} (E\ot \ov{H}^s\ot \ov{A}^r\ot k) \sub E\ot\ov{E}^{n+1}\ot k$,
we have that, on $k\ot \ov{H}^s\ot \ov{A}^{n+1-s}\ot k$,
$$
\align
\psi_{n+1}\circ \phi_{n+1} & = \ov{\si}_{n+1}\circ\psi_n\circ b'_{n+1}\circ
\phi_{n+1}\\
&= \ov{\si}_{n+1}\circ\psi_n\circ b'_{n+1}\circ \xi_{n+1}\circ\phi_n\circ
d_{n+1}\\
& = \ov{\si}_{n+1}\circ\psi_n\circ\phi_n\circ d_{n+1} -
\ov{\si}_{n+1}\circ\psi_n\circ \xi_n\circ b'_n\circ\phi_n \circ d_{n+1} \\
& = \ov{\si}_{n+1}\circ d_{n+1} =  id_{n+2} - d_{n+2}\circ  \ov{\si}_{n+2}.
\endalign
$$
So, to finish the proof it suffices to check that $\ov{\si}_{n+2}(k\ot
\ov{H}^s\ot \ov{A}^{n+1-s}\ot k) = 0$, which follows easily from the
definition of $\ov{\si}_*$\qed

\enddemo

Let $F^i(X_n) = \bigoplus_{0\le s\le i} E\ot \ov{H}^s \ot \ov{A}^{n-s}\ot
E$ and let $F^i(\B_n(E))$ be the sub-bimodule of $\B_n(E)$ generated by the
tensors $1\ot x_1\ot\cdots\ot x_n\ot 1$ such that at least $n-i$ of the
$x_j$'s belong to $A$.  The normalized Hochschild resolution $(\B_*(E)
,b'_*)$ and the resolution $(X_*, d_*)$ are filtered by $F^0(\B_*(E)) \sub
F^1(\B_*(E)) \sub F^2(\B_*(E)) \sub\dots$ and $F^0(X_*) \sub F^1(X_*) \sub
F^2(X_*) \sub\dots$, respectively

\proclaim{Proposition 1.2.2} The maps $\phi_*$, $\psi_*$ and $\om_{*+1}$
preserve filtrations.
\endproclaim

\demo{Proof} Let $Q^i_j = E\ot \ov{H}^i\ot \ov{A}^{n-j}\ot k$. We claim
that

\smallskip

\item{a)} $\ov{\si}_{n+1}(F^i(X_n)) \sub F^i(X_{n+1})$ for all $0\le i<n$,

\smallskip

\item{b)} $\ov{\si}_{n+1}(E\ot \ov{H}^i\ot \ov{A}^{n-i}\ot A) \sub
Q^i_{i-1} + F^{i-1}(X_{n+1})$ for all $0\le i\le n$,

\smallskip

\item{c)} $\ov{\si}_{n+1}(E\ot \ov{H}^n\ot E) \sub E\ot \ov{H}^{n+1}\ot
k + F^n(X_{n+1})$ for all $n\ge 0$,

\smallskip

\item{d)} $\psi_n(F^i(\B_n(E))\cap E\ot \ov{E}^n\ot k) \sub
Q^i_i + F^{i-1}(X_n)$.

\smallskip

In fact a), b) and c) follow immediately from the definition of
$\ov{\si}_{n+1}$. Suppose d) is valid for $n$. Let $\bx =
x_0\ot\cdots\ot x_{n+1}\ot 1 \in F^i(\B_{n+1}(E))\cap E\ot \ov{E}^{n+1}\ot
k$. Using a) and b), we get that for $1\le j \le n$,
$$
\ov{\si}_{n+1}(\psi_n(\bx_{0,j-1}\ot x_jx_{j+1}\ot \bx_{j+2,n+1} \ot 1))
\sub \ov{\si}_{n+1} (Q^i_i\! + F^{i-1}\!(X_n))
 \sub Q^i_{i-1}\!+ F^{i-1}\!(X_n).
$$
Since $\psi_{n+1}(\bx) = \ov{\si}_{n+1}\circ\psi_n \circ b'_{n+1}
(\bx)$, to prove d) for $n+1$ we only must check that $\ov{\si}_{n+1}
(\psi_n(\bx_{0,n+1})) \sub Q^i_{i-1} + F^{i-1}(X_n)$. If $x_{n+1}\in A$, then using
a) and b), we get
$$
\align
\ov{\si}_{n+1}(\psi_n(\bx_{0,n+1})) & = \ov{\si}_{n+1}(\psi_n(\bx_{0n}\ot
1_E)x_{n+1})\\
&\sub \ov{\si}_{n+1}(E\ot\ov{H}^i \ot \ov{A}^{n-i}\ot A + F^{i-1}(X_n))\\
&\sub Q^i_{i-1} + F^{i-1}(X_n),
\endalign
$$
and if $x_{n+1}\notin A$, then $\bx_{0,n+1} \in F^{i-1}(\B_n(E))$, which
together a) and c), implies that
$$
\ov{\si}_{n+1}(\psi_n(\bx_{0,n+1})) \sub \ov{\si}_{n+1}(F^{i-1}(X_n)) \sub
Q^i_{i-1} + F^{i-1}(X_{n+1}).
$$
From d) follows immediately that $\psi_*$ preserves filtrations. Next,
assuming that $\phi_n$ preserve filtrations, we prove that $\phi_{n+1}$
does it. Let $\bx\in F^i(X_{n+1})\cap Q^i_{i-1}$. Since $\phi_{n+1}(\bx) =
\xi_{n+1}\circ \phi_n\circ d_n(\bx)$ and
$$
\xi_{n+1}(\phi_n(d_{rs}^l(\bx))) \sub\xi_{n+1}(\phi_n(F^{i-l}\!(X_n)))
\sub\xi_{n+1}(F^{i-l}\!(\B_{@!n}(E))) \sub F^{i-l+1}\!(\B_{@!n+1}(E)),
$$
it suffices to see that $\xi_{n+1}(\phi_n(d_{rs}^0(\bx))) \sub
F^i(\B_{n+1}(E))$ for $\bx=1\ot \bh\ot \ba \ot 1$, with $\bh=
h_1\ot\cdots\ot h_i$ and $\ba= a_1\ot \cdots\ot a_{n+1-i}$. Since
$\phi_n(Q^i_i) \sub E\ot \ov{E}^n \ot k$, we have
$$
\align
\xi_{n+1} \circ \phi_n\circ d_{rs}^0 (\bx)
&= (-1)^r\xi_{n+1}\circ\phi_n (1\ot \bh\ot \ba)\\
& = (-1)^r\xi_{n+1}(\phi_n(1\ot \bh\ot \ba_{1,n-i}\ot 1)a_{n+1-i})\\
& \sub \xi_{n+1}(F^i(\B_n(E))\cap E\ot\ov{E}^n\ot A)\\
& \sub F^i(\B_{n+1}(E)).
\endalign
$$
Next, we prove that $\om_*$ preserves filtrations. Assume that $\om_n$ does
it. Let $\bx = x_0\ot\cdots\ot x_n\ot 1\in F^i(\B_n(E)) \cap E\ot
\ov{E}^n\ot k$. It is evident that $\om_{n+1}(\bx) = \xi_{n+1}\circ
\phi_n\circ \psi_n(\bx) - \xi_{n+1}\circ\om_n\circ b'_n(\bx)$. Since
$\xi_{n+1}(\phi_n(Q^i_i)) \sub \xi_{n+1}(E\ot \ov{E}^n \ot k) = 0$, from d)
we get
$$
\xi_{n+1}\circ\phi_n\circ \psi_n(\bx)  \in \xi_{n+1}\circ \phi_n\bigl(Q^i_i
+ F^{i-1}(X_n)\bigr) \sub\xi_{n+1}(F^{i-1}(\B_n(E)))\sub F^i(\B_n(E)).
$$
It remains to check that $\xi_{n+1}\circ\om_n\circ b'_n(\bx) \sub
F^i(\B_n(E))$. Since $\om_n(E\ot \ov{E}^{n-1}\ot k) \sub E\ot\ov{E}^n\ot
k$, we have $\xi_{n+1}\circ\om_n \circ b'_n(\bx) =
(-1)^{n-1}\xi_{n+1}\circ\om_n (\bx_{0n})$. Hence, if $x_n\in A$, then
$$
\align
\xi_{n+1}\circ\om_n\circ b'_n(\bx) &= (-1)^{n-1}\xi_{n+1}(
\om_n(\bx_{0,n-1}\ot 1)x_n)\\
&\sub \xi_{n+1}(F^i(\B_n(E))\cap E \ot \ov{E}^n\ot A)\\
&\sub F^i(\B_{n+1}(E)),
\endalign
$$
and if $x_n\notin A$, then $\bx\in F^{i-1}(\B_{n-1}(E))$, and so
$$
\xi_{n+1}\circ\om_n\circ b'_n(\bx)  = (-1)^{n-1} \xi_{n+1}
\circ\om_n(\bx_{0n}) \sub\xi_{n+1}(F^{i-1}(\B_n(E))) \sub F^i(\B_{n+1}
E)\qed
$$
\enddemo

\head 2. The Hochschild (co)homology of a crossed product \endhead
Let $E=A\#_f H$ and $M$ an $E$-bimodule. In this section we use
Theorem~1.1.1 in order to construct complexes $\wh{X}_*(E,M)$ and
$\wh{X}_*(E,M)$, simpler than the canonical ones, giving the Hochschild
homology and cohomology of $A$ with coefficients in $M$ respectively. These
complexes have natural filtrations that allow us to obtain spectral
sequences converging to $\H_*(E,M)$ and $\H^*(E,M)$ respectively.

\specialhead 2.1. Hochschild homology \endspecialhead
Let $\wh{d}_{rs}^l\: M\ot \ov{H}^s \ot \ov{A}^r \to M\ot \ov{H}^{s-l} \ot
\ov{A}^{r+l-1}$ ($r,s\ge 0$, $0\le l \le s$ and $r+l>0$) be the morphisms
defined by:
$$
\align
&\wh{d}_{rs}^0(\bx) = ma_1^{\ov{\bh^{(1)}}}
\ot \bh^{(2)}\ot \ba_{2r} + (-1)^r a_rm\ot \bh\ot \ba_{1,r-1}\\
&\phantom{\wh{d}_{rs}^0(\bx)} + \sum_{i=1}^{r-1} (-1)^i
m\ot \bh\ot \ba_{1,i-1}\ot a_ia_{i+1}\ot \ba_{i+2,r},\\
\vspace{1.5\jot}
&\wh{d}_{rs}^1(\bx) = (-1)^r m(1\# h_1)\ot \bh_{2s}\ot \ba +
(-1)^{r+s}(1\#h_s^{(2)}) m \ot \bh_{1,s-1}^{(1)}\ot \ba^{h_s^{(1)}}\\
&\phantom{\wh{d}_{rs}^1(\bx)} + \sum_{i=1}^{s-1} (-1)^{r+i}
mf(h_i^{(1)},h_{i+1}^{(1)})^{\bh_{1,i-1}^{(1)}} \ot \bh_{1,i-1}^{(2)}\ot
h_i^{(2)}h_{i+1}^{(2)} \ot \bh_{i+2,s}\ot \ba\\
\vspace{1.5\jot}
&\wh{d}_{rs}^l(\bx) = (-1)^{l(r+s)}(1\#\fh_{s-l+1,s}^{(2)})
m\ot\bh_{1,s-l}\ot F_r^{(l)}(\bh_{s-l+1,s}^{(1)}\ot \ba),
\endalign
$$
where $\bx=m\ot\bh\ot\ba$, with $\ba = a_1\ot\cdots\ot a_r$ and $\bh =
h_1\ot\cdots\ot h_s$.

\proclaim{Theorem 2.1.1} The Hochschild homology of $E$ with coefficients
in $M$ is the homology of the chain complex
$$
\wh{X}_*(E,M) = \quad \wh{X}_0 @<\wh{d}_1 << \wh{X}_1 @<\wh{d}_2 <<
\wh{X}_2 @<\wh{d}_3 << \wh{X}_3 @<\wh{d}_4 << \wh{X}_4 @<\wh{d}_5 <<
\wh{X}_5 @<\wh{d}_6 << \wh{X}_6 @<\wh{d}_7 << \dots,
$$
where $\displaystyle{\wh{X}_n = \bigoplus_{r+s=n} M\ot\ov{H}^s\ot
\ov{A}^r}$ and $\displaystyle{\wh{d}_n = \sum_{r+s=n\atop r+l> 0}
\sum^s_{l=0} \wh{d}^l_{rs}}$.
\endproclaim

\demo{Proof} It follows from the fact that $\wh{X}_*(E,M)\simeq
M\ot_{E^e}(X_*,d_*)$. An isomorphism is provided by the maps
$\wh{\theta}_{rs}\: M\ot\ov{H}^s \ot \ov{A}^r @>>> M\ot_{E^e} X_{rs} $,
defined by $\wh{\theta}_{rs}(m\ot \bh\ot \ba) = m \ot (1_E\ot \bh \ot
\ba\ot 1_E)$\qed
\enddemo

\subhead 2.1.2. A spectral sequence \endsubhead Let $F^i(\wh{X}_n) =
\bigoplus_{0\le s\le i} M\ot \ov{H}^s \ot \ov{A}^{n-s}$. Clearly
$F^0(\wh{X}_*) \sub F^1(\wh{X}_*) \sub\dots$ is a filtration of
$\wh{X}_*(E,M)$. Using this fact we obtain:

\proclaim{Corollary 2.1.2.1} There is a convergent spectral sequence
$$
E^1_{rs} = \H_r(A,M\ot \ov{H}^s) \Rightarrow \H_{r+s}(E,M),
$$
where $M\ot \ov{H}^s$ is considered as an $A$-bimodule via $a_1(m \ot
\bh_{1s})a_2 = a_1 ma_2^{\ov{\bh_{1s}^{(1)}}} \ot \bh_{1s}^{(2)}$.
\endproclaim

The normalized Hochschild complex $(M\ot \ov{E}^*,b_*)$ has a filtration
$F^0(M\ot \ov{E}^*) \sub F^1(M\ot \ov{E}^*) \sub F^2(M\ot \ov{E}^*)
\sub\dots$, where $F^i(M\ot \ov{E}^n)$ is the $k$-submodule of $M\ot
\ov{E}^n$ generated by the tensors $m\ot x_1\ot\cdots\ot x_n$ such that at
least $n-i$ of the $x_j$'s belong to $A$. The spectral sequence associate
to this filtration is called the homological Hochschild-Serre spectral
sequence. Since, for each extension of groups $N\sub G$ with $N$ a normal
subgroup, it is hold that $k[G]$ is a crossed product of $k[G/N]$ on
$k[N]$, the following theorem (joint with Corollary~3.1.3 below) gives, as
a particular case, the homological version of the main results of
\cite{H-S}.

\proclaim{Theorem 2.1.2.2} The homological Hochschild-Serre spectral
sequence is isomorphic to the one obtained in Corollary~2.1.2.1.
\endproclaim

\demo{Proof} It is an easy consequence of Propositions~1.2.1 and 1.2.2.
\enddemo

\subhead 2.1.3. A decomposition of $\wh{X}_*(E,M)$ \endsubhead Let $[H,H]$
be the $k$-submodule of $H$ spanned by the set of all elements $ab-ba$
($a,b\in H$). It is easy to see that $[H,H]$ is a coideal in $H$. Let
$\breve H$ be the quotient coalgebra $H/[H,H]$. Given $h\in H$, we let
$[h]$ denote the class of $h$ in $\breve H$. Given a subcoalgebra $C$ of
$\breve H$ and a right $\breve H$-comodule $(N,\rho)$, we put $N^C = \{n\in
N\:\rho(n)\in N\ot C\}$. It is well known that if $\breve H$ decomposes as
a direct sum of subcoalgebras $C_i$ ($i\in I$), then $N= \bigoplus_{i\in I}
N^{C_i}$.

Now, let us assume that $M$ is a Hopf bimodule. That is, $M$ is an
$E$-bimodule and a right $H$-comodule, and the coaction $m\mapsto
m^{(0)}\ot m^{(1)}$ verifies:
$$
((a\# h)m(b\# l))^{(0)}\! \ot ((a\# h)m(b\# l))^{(1)} = (a\# h^{(1)}\!)
m^{(0)}\! (b\# l^{(1)}\!) \ot (a\# h^{(2)}\!) m^{(1)}\! (b\# l^{(2)}\!).
$$
For each $n\ge 0$, $\wh{X}_n$ is an $\breve H$-comodule via
$$
\rho_{n}(m\ot\bh_{1s}\ot\ba_{1r}) = m^{(0)}\ot\bh_{1s}^{(1)}\ot\ba_{1r} \ot
[m^{(1)}\fh_{1s}^{(2)}] \qquad\text{($r+s=n$).}
$$
Moreover, the map $\rho_*\:\wh{X}_*(E,M) @>>> \wh{X}_*(E,M)\ot \breve H$ is
a map of complexes. This fact implies that if $C$ is a subcoalgebra of
$\breve H$, then $\wh{d}_n(\wh{X}_n^C)\sub \wh{X}_{n-1}^C$. We consider
the subcomplex $\wh{X}_*^C(E,M)$ of $\wh{X}_*(E,M)$, with modules
$\wh{X}_n^C$, and we let $\H_*^C(E,M)$ denote its homology. By the
above discussion, if $\breve H$ decomposes as a direct sum of subcoalgebras
$C_i$ ($i\in I$), then $\wh{X}_*(E,M)= \bigoplus_{i\in I} \wh{X}_*^{C_i}
(E,M)$. Consequently $\H_*(E,M)= \bigoplus_{i\in I} \H_*^{C_i}(E,M)$.
Finally, the filtration of $\wh{X}_*(E,M)$ induces a filtration on
$\wh{X}_*^C(E,M)$. Hence we have a convergent spectral sequence
$$
E^1_{rs} = \H_r(A,(M\ot \ov{H}^s)^C) \Rightarrow \H_{r+s}^C(E,M),
$$
where $(M\ot \ov{H}^s)^C$ is an $A$-bimodule via $a_1(m\ot \bh_{1s})a_2 =
a_1 ma_2^{\ov{\bh_{1s}^{(1)}}} \ot \bh_{1s}^{(2)}$.

\subhead 2.1.4. Compatibility with the canonical decomposition \endsubhead
Let us assume that $k\supseteq \bQ$, $H$ is cocommutative, $A$ is
commutative, $M$ is symmetric as an $A$-bimodule and the cocycle $f$ takes
its values in $k$. In \cite{G-S1} was obtained a decomposition of the
canonical Hochschild complex $(M\ot \ov{A}^*,b_*)$. It is easy to check
that the maps $\wh{d}_0$ and $\wh{d}_1$ are compatible with this
decomposition. Since $\wh{d}_l = 0$ for all $l\ge 2$, we obtain a
decomposition of $\wh{X}_*(E,M)$, and then a decomposition of $\H_*(E,M)$.

\specialhead 2.2. Hochschild cohomology \endspecialhead
Let $\wh{d}^{rs}_l\:\Hom_k(\ov{H}^{s-l}\ot \ov{A}^{s+l-1},M)\to
\Hom_k(\ov{H}^s \ot \ov{A}^r,M)$ ($0\le l \le s$, $r+l>0$) be the morphisms
defined by:
$$
\align
&\wh{d}^{rs}_0(\varphi)(\bx) = a_1^{\bh^{(1)}}
\varphi(\bh^{(2)}\ot \ba_{2r}) + (-1)^r \varphi(\bh\ot
\ba_{1,r-1})a_r\\
&\phantom{\wh{d}^{rs}_0(\varphi)(\bx)} + \sum_{i=1}^{r-1}
(-1)^i \varphi(\bh\ot \ba_{1,i-1}\ot a_ia_{i+1}\ot \ba_{i+2,r}),\\
\vspace{1.5\jot}
&\wh{d}^{rs}_1(\varphi)(\bx) = (-1)^r (1\# h_1)
\varphi(\bh_{2s}\ot\ba) + (-1)^{r+s} \varphi(\bh_{1,s-1}\ot
\ba^{h_s^{(1)}})(1\# h_s^{(2)})\\
&\phantom{\wh{d}^{rs}_1(\varphi)(\bx)} + \sum_{i=1}^{s-1}
(-1)^{r+i} f(h_i^{(1)},h_{i+1}^{(1)})^{\bh_{1,i-1}^{(1)}} \varphi(
\bh_{1,i-1}^{(2)}\ot h_i^{(2)}h_{i+1}^{(2)}\ot \bh_{i+2,s}\ot\ba),\\
\vspace{1.5\jot}
&\wh{d}^{rs}_l(\varphi)(\bx) = (-1)^{l(r+s)}\varphi \bigl( \bh_{1,s-l} \ot
F_r^{(l)}(\bh_{s-l+1,s}^{(1)}\ot \ba)\bigr)(1\# \fh_{s-l+1,s}^{(2)}),
\endalign
$$
where $\bx=\bh\ot\ba$, with $\ba = a_1\ot\cdots\ot a_r$ and $\bh =
h_1\ot\cdots\ot h_s$.

\proclaim{Theorem 2.2.1} The Hochschild cohomology of $E$ with coefficients
in $M$ is the homology of
$$
\wh{X}^*(E,M)=\quad \wh{X}^0 @>\wh{d}^1>> \wh{X}^1 @>\wh{d}^2>> \wh{X}^2
@>\wh{d}^3>> \wh{X}^3 @>\wh{d}^4>> \wh{X}^4 @>\wh{d}_5>> \wh{X}^5
@>\wh{d}^6>> \wh{X}^6 @>\wh{d}^7>> \dots,
$$
where $\displaystyle{\wh{X}^n=\bigoplus_{r+s=n}\Hom_k(\ov{H}^s\ot \ov{A}^r,
M)}$ and $\displaystyle{\wh{d}^n = \sum_{r+s=n\atop r+l>0} \sum_{l=0}^s
\wh{d}_l^{rs}}$.
\endproclaim

\demo{Proof} It follows from the fact that $\wh{X}^*(E,M)\simeq
\Hom_{E^e}((X_*,d_*),M)$. An isomorphism is provided by the maps
$\wh{\theta}^{rs}\: \Hom_k(\ov{H}^s\ot \ov{A}^r,M) @>>> \Hom_{E^e}
(X_{rs},M)$, defined by $\wh{\theta}^{rs}(\varphi)(1_E\ot\bx \ot 1_E) =
\varphi(\bx)$\qed
\enddemo

\subhead 2.2.2. A spectral sequence \endsubhead Let ${F^i(\wh{X}^n) =
\bigoplus_{s\ge i} \Hom_k(\ov{H}^s \ot \ov{A}^{n-s},M)}$. Clear\-ly $F_0
\sup F_1 \sup F_2 \sup\dots$  is a filtration of $\wh{X}^*(E,M)$. Using
this fact we obtain:

\proclaim{Corollary 2.2.2.1} There is a convergent spectral sequence
$$
E_1^{rs} = \H^r(A,\Hom_k(\ov{H}^s,M)) \Rightarrow \H^{r+s}(E,M),
$$
where $\Hom_k(\ov{H}^s,M)$ is considered as an $A$-bimodule via $a_1
\varphi a_2(\bh) = a_1^{\bh^{(1)}} \varphi(\bh^{(2)}) a_2$.
\endproclaim

Let $F_i(\Hom_k(\ov{E}^n,M))$ be the $k$-submodule of $(\Hom_k(\ov{E}^*,M),
b^*)$ consisting of maps $f\in \Hom_k(\ov{E}^n,M)$, for which $f(x_1\ot
\cdots\ot x_n) = 0$ whenever $n-i$ of the $x_j$'s belong to $A$. The
normalized Hochschild complex $(\Hom_k(\ov{E}^*,M),b^*)$ is filtered by
$F_0(\Hom_k (\ov{E}^*, M)) \sup F_1(\Hom_k(\ov{E}^*,M)) \sup
F_2(\Hom_k(\ov{E}^*,M)) \sup \dots$. The spectral sequence associated to
this filtration is called the cohomological  Hochschild-Serre spectral
sequence. The following theorem (joint with Corollary~3.2.3 below) gives,
as a particular case, of the main results of \cite{H-S}.

\proclaim{Theorem 2.2.2.2} The cohomological Hochschild-Serre spectral
sequence is isomorphic to the one obtained in Corollary~2.2.2.1.
\endproclaim

\demo{Proof} It is an easy consequence of Propositions~1.2.1 and 1.2.2.
\enddemo

\subhead 2.2.3. Compatibility with the canonical decomposition \endsubhead
Assume that $k\supseteq \bQ$, $H$ is cocommutative, $A$ is commutative, $M$
is symmetric as an $A$-bimodule and the cocycle $f$ takes its values in
$k$. Then, the Hochschild cohomology $\H^*(E,M)$ has a decomposition
similar to the one obtained in 2.1.4 for the Hochschild homology.

\head 3. The Hochschild (co)homology of a crossed product with
invertible cocycle \endhead
Let $E=A\#_f H$ and $M$ an $E$-bimodule. Assume that the cocycle $f$ is
invertible. Then, the map $h\mapsto 1\# h$ is convolution invertible and
its inverse is the map $h\mapsto (1\# h)^{-1} = f^{-1}(S(h^{(2)}),h^{(3)})
\# S(h^{(1)})$. Under this hypothesis, we prove that the complexes
$\wh{X}_*(E,M)$ and $\wh{X}^*(E,M)$ of Section~2 are isomorphic to simpler
complexes. These complexes have natural filtrations, which give the
spectral sequences obtained in \cite{S}. Using these facts and a theorem of
Gerstenhaber and Schack, we prove that if the $2$-cocycle $f$ takes its
values in a separable subalgebra of $A$, then the Hochschild (co)homology
of $E$ with coefficients in $M$ is the (co)homology of $H$ with
coefficients in a (co)chain complex. Finally, as an application we obtain
some results about the $\Tor_*^E$ and $\Ext^*_E$ functors and an upper
bound for the global dimension of $E$.

\specialhead 3.1. Hochschild homology \endspecialhead
Let $\ov{d}_{rs}^l\: M\ot \ov{A}^r \ot \ov{H}^s \to M\ot \ov{A}^{r+l-1} \ot
\ov{H}^{s-l}$ ($r,s\ge 0$, $0\le l \le s$ and $r+l>0$) be the morphisms
defined by:
$$
\allowdisplaybreaks
\align
&\ov{d}_{rs}^0(\bx) = ma_1\ot \ba_{2r}\ot \bh + (-1)^r a_rm\ot
\ba_{1,r-1}\ot \bh\\
&\phantom{\ov{d}_{rs}^0(\bx)} + \sum_{i=1}^{r-1} (-1)^i
m\ot \ba_{1,i-1}\ot a_ia_{i+1}\ot \ba_{i+2,r}\ot \bh,\\
\vspace{1.5\jot}
&\ov{d}_{rs}^1(\bx)= (-1)^r m\ep(h_1)\ot \ba\ot\bh_{2s} +
(-1)^{r+s}(1\#h_s^{(3)}) m (1\# h_s^{(1)})^{-1}\ot\ba^{h_s^{(2)}}\ot
\bh_{1,s-1} \\
&\phantom{\ov{d}_{rs}^1(\bx)} + \sum_{i=1}^{s-1} (-1)^{r+i} m \ot
\ba\ot \bh_{1,i-1}\ot h_ih_{i+1} \ot \bh_{i+2,s},\\
\vspace{1.5\jot}
&\ov{d}_{rs}^l(\bx) = (-1)^{l(r+s)}(1\#\fh_{s-l+1,s}^{(3)}) m
(1\#\fh_{s-l+1,s}^{(1)})^{-1}\ot F_r^{(l)}(\bh_{s-l+1,s}^{(2)}\ot \ba),
\endalign
$$
where $\bx=m\ot \ba\ot\bh$, with $\ba = a_1\ot\cdots\ot a_r$ and $\bh=
h_1\ot\cdots\ot h_s$. Let $\ov{X}_*(E,M)$ be the complex
$$
\ov{X}_*(E,M)=\quad \ov{X}_0 @<\ov{d}_1 << \ov{X}_1 @<\ov{d}_2 <<
\ov{X}_2 @<\ov{d}_3 << \ov{X}_3 @<\ov{d}_4 << \ov{X}_4 @<\ov{d}_5 <<
\ov{X}_5 @<\ov{d}_6 << \ov{X}_6 @<\ov{d}_7 << \dots,
$$
where
$\ov{X}_n = \bigoplus_{r+s=n} M\ot\ov{A}^r\ot\ov{H}^s$ and $\ov{d}_n =
\sum_{r+s=n\atop r+l> 0} \sum^s_{l=0} \ov{d}^l_{rs}$.

\proclaim{Theorem 3.1.1} The map $\theta_*\:\wh{X}_*(E,M) @>>>
\ov{X}_*(E,M)$, given by
$$
\theta_n(m\ot \bh\ot \ba) = m (1\# h_1^{(1)}) \cdots (1\#
h_s^{(1)}) \ot \ba \ot \bh^{(2)}\qquad \text{($r+s = n$)},
$$
is an isomorphism of complexes. Consequently, the Hochschild homology of
$E$ with coefficients in $M$ is the homology of $\ov{X}_*(E,M)$.

\endproclaim

\demo{Proof} A direct computation shows that $\theta_*$ is a morphism of
complexes. The inverse map of $\theta_n$ is the map $m\ot \ba\ot
\bh \mapsto  m(1\# h_s^{(1)})^{-1}\cdots (1\# h_1^{(1)})^{-1}\ot
\bh^{(2)}\ot \ba$\qed
\enddemo

Note that when $f$ takes its values in $k$, then $\ov{X}_*(E,M)$ is the
total complex of the double complex $\bigl(M\ot \ov{A}^* \ot \ov{H}^*,
\ov{d}_{**}^0,\ov{d}_{**}^1\bigr)$.

\medskip

For each $h\in H$, we have the morphism $\theta^h_*\: (M\ot\ov{A}^*,b_*)
@>>> (M\ot\ov{A}^*,b_*)$, defined by $\theta^h_r(m\ot \ba) =
(1\#h_s^{(3)}) m (1\# h_s^{(1)})^{-1}\ot\ba^{h_s^{(2)}}$.

\proclaim{Proposition 3.1.2} For each $h,l\in H$ the endomorphisms of
$\H_*(A,M)$ induced by $\theta^h_*\circ \theta^l_*$ and by $\theta^{hl}_*$
coincide. Consequently $\H_*(A,M)$ is a left $H$-module.
\endproclaim

\demo{Proof} By a standard argument it is sufficient to prove it for
$\H_0(A,M)$, and in this case the result is immediate\qed
\enddemo

\proclaim{Corollary 3.1.3} The chain complex $\ov{X}_*(E,M)$ has a
filtration $F^0 \sub F^1 \sub\dots$, where $F^i(\ov{X}_n) = \bigoplus_{
0\le s\le i} M\ot \ov{A}^{n-s} \ot \ov{H}^s$. The spectral sequence of this
filtration is isomorphic to the one obtained in Corollary~2.1.2. From
Proposition~3.1.2 it follows that if $H$ is a flat $k$-module, then
$E^1_{rs} = \H_r(A,M)\ot \ov{H}^s$ and $E^2_{rs} = \H_s(H,\H_r(A,M))$.

\endproclaim

Given an $A$-bimodule $M$ we let $[A,M]$ denote the $k$-submodule of $M$
generated by the commutators $am - ma$ ($a\in A$ and $m\in M$).

\remark{Remark 3.1.4} From Corollary~3.1.3 it follows immediately that if
$A$ is separable, then $\H_*(E,M) = \H_*(H,M/[A,M])$, and if $A$ is
quasi-free, then there is a long exact sequence
$$
\align
\dots @>>> &\H_{n+1}(H,\H_0(A,M)) @>>>\H_{n-1}(H,\H_1(A,M)) @>>> \H_n(E,M)
@>>>\\
&\H_n(H,\H_0(A,M)) @>>> \H_{n-2}(H,\H_1(A,M)) @>>> \H_{n-1} (E,M) @>>>\dots.
\endalign
$$
\endremark

\subhead 3.1.5. Separable subalgebras \endsubhead Let $S$ be a separable
subalgebra of $A$. Next we prove that if the $2$-cocycle $f$ takes its
values in $S$, then the Hochschild homology of $E$ with coefficients in $M$
is the homology of $H$ with coefficients in a chain complex. When $S$
equals $A$ we recover the first part of Remark~3.1.4. Assume that
$f(h,l)\in S$ for all $h,l\in H$. Let $\wt{A} = A/S$, $\wt{A}^0 = S$ and
$\wt{A}^r = \wt{A}\ot_S\cdots\ot_S \wt{A}$ ($r$-times) for $r>0$, and let
$M\ot_S \wt{A}^r\ot_S = M\ot_{A^e} (A\ot_s\wt{A}^r \ot_s A) = M\ot_S
\wt{A}^r\ot_{S^e} S$ be the cyclic tensor product over $S$ of $M$ and
$\wt{A}^r$ (see \cite{G-S2} or \cite{Q}). Using the fact that $f$ takes its
values in $S$, it is easy to see that $H$ acts on $(M\ot_S \wt{A}^r
\ot_S,b_*)$ via
$$
h\cdot \bigl(m\ot_S \wt{\ba}\bigr) = (1\#h^{(3)}) m (1\# h^{(1)})^{-1}\ot_S
\wt{\ba}^{\ov{h^{(2)}}},
$$
where $m\ot_S \wt{\ba} = m\ot_S a_1 \ot_S \cdots \ot_S a_r \ot_S$ and
$\wt{\ba}^{\ov{h^{(2)}}} = a_1^{\ov{h^{(2)}}} \ot_S \cdots \ot_S
a_r^{\ov{h^{(r+1)}}} \ot_S$.

\proclaim{Theorem 3.1.5.1} The Hochschild homology $\H_*(E,M)$, of $E$ with
coefficients in $M$, is the homology of $H$ with coefficients in $(M\ot_S
\wt{A}^r \ot_S,b_*)$.

\endproclaim

\demo{Proof} Let $((M\ot_S \wt{A}^* \ot_S) \ot \ov{H}^*,\wt{d}^0_{**},
\wt{d}^1_{**})$ be the double complex with horizontal differentials
$$
\align
&\wt{d}_{rs}^0(\bx) = ma_1\ot_S \wt{\ba}_{2r}\ot \bh + (-1)^r a_rm\ot_S
\wt{\ba}_{1,r-1}\ot \bh\\
&\phantom{\wt{d}_{rs}^0(\bx)} + \sum_{i=1}^{r-1} (-1)^i m\ot_S
\wt{\ba}_{1,i-1}\ot_S a_ia_{i+1}\ot_S\wt{\ba}_{i+2,r}\ot \bh,\\
\intertext{and vertical differentials}
&\wt{d}_{rs}^1(\bx) = (-1)^r m\ot_s \wt{\ba} \ot \bh_{2s} +(-1)^{r+s} (1\#
h_s^{(3)}) m (1\# h_s^{(1)})^{-1} \ot_S\wt{\ba}^{h_s^{(2)}} \ot \bh_{1,s-1}
\\
&\phantom{\wt{d}_{rs}^1(\bx)} + \sum_{i=1}^{s-1} (-1)^{r+i} m\ot_S\wt{\ba}
\ot \bh_{1,i-1}^{(2)}\ot h_i^{(2)} h_{i+1}^{(2)} \ot \bh_{i+2,s},\\
\endalign
$$
where $\bx = m\ot \wt{\ba}\ot \bh$, with $\wt{\ba} = a_1\ot_S\cdots\ot_S
a_r\ot_S$ and $\bh = h_1\ot\cdots\ot h_s$. Let $\ov{X}^S_*(E,M)$ be the
total complex of $((M\ot_S \wt{A}^* \ot_S) \ot \ov{H}^*,\wt{d}^0_{**},
\wt{d}^1_{**})$. We must prove that $\H_*(E,M)$ is the homology of
$\ov{X}^S_*(E,M)$. Let $\pi_*\: \ov{X}_*(E,M)\to \ov{X}^S_*(E,M)$ be the
map $m\ot \ba_{1r}\ot \bh_{1s}\mapsto m\ot_S \wt{\ba}_{1r}\ot \bh_{1s}$.
Consider the filtration $F^{0S}_* \sub F^{1S}_* \sub F^{2S}_* \sub\dots$ of
$\ov{X}^S_*(E,M)$, where $F^{iS}_n  = \bigoplus_{0\le s\le i} (M\ot_S
\wt{A}^{n-s}\ot_S)\ot \ov{H}^s$. From Theorem~1.2 of \cite{G-S2}, it follows
that $\pi_*$ is a morphism of filtered complexes inducing an
quasi-isomorphism between the graded complexes associated to the
filtrations of $\ov{X}_*(E,M)$ and $\ov{X}^S_*(E,M)$. Consequently $\pi_*$
is a quasi-isomorphism. The proof can be finished by applying
Theorem~3.1.1\qed

\enddemo

\subhead 3.1.6. A decomposition of $\ov{X}_*(E,M)$ \endsubhead Here we
freely use the notations of Subsection~2.1.3. Suppose $M$ is a Hopf
bimodule. A direct computation shows that the $\breve H$-coaction of
$\ov{X}_*(E,M)$, obtained transporting the one of $\wh{X}_*(E,M)$ through
$\theta_*\:\wh{X}_*(E,M) @>>> \ov{X}_*(E,M)$, is given by
$$
m\ot \ba\ot \bh \mapsto m^{(0)}\ot \ba\ot \bh^{(2)} \ot
m^{(1)}S(h_s^{(1)})\cdots S(h_1^{(1)})h_1^{(3)}\cdots h_s^{(3)}.\tag 2
$$
For each subcoalgebra $C$ of $\breve H$, we consider the subcomplex
$\ov{X}_*^C(E,M)$ of $\ov{X}_*(E,M)$ with modules $\ov{X}^C_n$, and we let
$\H_*^C(E,M)$ denote its homology. If $\breve H$ decomposes as a direct sum
of subcoalgebras $C_i$ ($i\in I$), then $\ov{X}_*(E,M) = \bigoplus_{i\in I}
\ov{X}_*^{C_i}(E,M)$. Consequently $\H_*(E,M)= \bigoplus_{i\in I}
\H_*^{C_i}(E,M)$. From $(2)$ it follows that if $\breve H$ is
cocommutative, then $\ov{X}_n^C = \bigoplus_{r+s=n} M^C\ot \ov{A}^r \ot
\ov{H}^s$. Finally, the filtration of $\ov{X}_*(E,M)$ induces a filtration
on $\ov{X}_*^C(E,M)$. Hence, when $\breve H$ is cocommutative and $H$ is a
flat $k$-module, we have a convergent spectral sequence
$$
E^2_{rs} = \H_r(H,\H_s(A,M^C)) \Rightarrow \H_{r+s}^C(E,M),
$$
where $\H_r(A,M^C)$ is a left $H$-module via the action introduced in
Proposition~3.1.2.

\subhead 3.1.7. An application to $\Tor_*^E$ \endsubhead Let $k$ be a
field, $B$ an arbitrary $k$-algebra, $M$ a right $B$-module and $N$ a left
$B$-module. It is well known that $\Tor_*^B(M,N)\simeq \H_*(B,N\ot M)$
(here $N\ot M$ is an $B$-bimodule via $a(n\ot m)b = an\ot mb$). This fact
and Corollary~3.1.3 show that if $k$ is a field, $M$ is a right $E$-module
and $N$ is a left $E$-module, then there is a convergent spectral sequence
$$
E^2_{rs} = \H_r(H,\Tor_s^A(M,N)) \Rightarrow \Tor_{r+s}^E(M,N).
$$

\specialhead 3.2. Hochschild cohomology\endspecialhead
Let $\ov{d}^{rs}_l\:\Hom_k(\ov{A}^{r+l-1}\ot \ov{H}^{s-l},M)\to
\Hom_k(\ov{A}^r \ot \ov{H}^s,M)$ ($0\le l \le s$, $r+l>0$) be the morphisms
defined by:
$$
\align
&\ov{d}^{rs}_0(\varphi)(\bx) = a_1\varphi(\ba_{2r}\ot \bh) + (-1)^r
\varphi(\ba_{1,r-1}\ot \bh)a_r\\
&\phantom{\ov{d}^{rs}_0(\varphi)(\bx)} + \sum_{i=1}^{r-1} (-1)^i
\varphi(\ba_{1,i-1}\ot a_ia_{i+1}\ot \ba_{i+2,r}\ot \bh),\\
\vspace{1.5\jot}
&\ov{d}^{rs}_1(\varphi)(\bx) = (-1)^r \ep(h_1)\varphi( \ba\ot \bh_{2s}) +
(-1)^{r+s} (1\# h_s^{(1)})^{-1} \varphi(\ba^{h_s^{(2)}}\ot\bh_{1,s-1})(1\#
h_s^{(3)}) \\
&\phantom{\ov{d}^{rs}_1(\varphi)(\bx)} + \sum_{i=1}^{s-1} (-1)^{r+i}
\varphi(\ba\ot\bh_{1,i-1}\ot h_ih_{i+1}\ot \bh_{i+2,s}),\\
\vspace{1.5\jot}
&\ov{d}^{rs}_l(\varphi)(\bx) = ((-1)^{l(r+s)}\# \fh_{s-l+1,s}^{(1)})^{-1}
\varphi\bigl(F_r^{(l)}(\bh_{s-l+1,s}^{(2)}\ot \ba) \ot \bh_{1,s-l} \bigr)
(1\# \fh_{s-l+1,s}^{(3)}),
\endalign
$$
where $\bx =\ba\ot \bh$, with $\ba= a_1\ot\cdots\ot a_r$ and
$\bh=h_1\ot\cdots\ot h_s$. Let $\ov{X}^*(E,M)$ be the complex
$$
\ov{X}^*(E,M)=\quad \ov{X}^0 @>\ov{d}^1>> \ov{X}^1 @>\ov{d}^2>> \ov{X}^2
@>\ov{d}^3>> \ov{X}^3 @>\ov{d}^4>> \ov{X}^4 @>\ov{d}_5>> \ov{X}^5
@>\ov{d}^6>> \ov{X}^6 @>\ov{d}^7 >> \dots,
$$
where $\ov{X}^n=\bigoplus_{r+s=n}\Hom_k(\ov{A}^r\ot\ov{H}^s,M)$ and
$\ov{d}^n = \sum_{r+s=n\atop r+l>0} \sum_{l=0}^s \ov{d}_l^{rs}$.

\proclaim{Theorem 3.2.1} The map $\theta^*\:\ov{X}^*(E,M) @>>>
\wh{X}^*(E,M)$, given by
$$
\theta^n(\varphi)(\bh\ot \ba) = (1\# h_1^{(1)})\cdots (1\#
h_s^{(1)}) \varphi(\ba\ot \bh^{(2)}) \qquad \text{($r+s = n$)},
$$
is an isomorphism of complexes. Consequently, the Hochschild cohomology of
$E$ with coefficients in $M$ is the homology of $\ov{X}^*(E,M)$.

\endproclaim

\demo{Proof} It is similar to the proof of Theorem~3.1.1\qed
\enddemo

Note that when $f$ takes its values in $k$, then $\ov{X}^*(E,M)$ is the
total complex of the double complex $\bigl(\Hom_k(\ov{A}^*\ot\ov{H}^*,M),
\ov{d}^{**}_0,\ov{d}^{**}_1\bigr)$.

\medskip

For each $h\in H$ we have the map $\theta_h^*\: (\Hom_k(\ov{A}^*,M), b^*)
@>>> (\Hom_k(\ov{A}^*,M),b^*)$, defined by $\theta_h^r(\varphi) (\ba)
= (1\# h^{(1)})^{-1} \varphi(\ba^{h^{(2)}})(1\# h^{(3)})$.

\proclaim{Proposition 3.2.2} For each $h,l\in H$ the endomorphisms of
$\H^*(A,M)$ induced by $\theta_l^*\circ \theta_h^*$ and by $\theta_{hl}^*$
coincide. Consequently $\H^*(A,M)$ is a right $H$-module.
\endproclaim

\demo{Proof} By a standard argument it is sufficient to prove it for
$\H^0(A,M)$, and in this case the result is immediate\qed
\enddemo

\proclaim{Corollary 3.2.3} The cochain complex $\ov{X}^*(E,M)$ has a
filtration $F_0 \sup F_1\sup \dots$, where $F_i(\ov{X}^n) = \bigoplus_{0
\le r\le n-i} \Hom_k(\ov{A}^r \ot \ov{H}^{n-r},M)$. The spectral sequence
of this filtration is isomorphic to the one obtained in Corollary~2.2.2.
From Proposition~3.2.2 it follows that $E_1^{rs} = \Hom_k(\ov{H}^s,
\H^r(A,M))$ and $E_2^{rs} = \H^s(H,\H^r(A,M))$.
\endproclaim

\medskip

Given an $A$-bimodule $M$, we let $M^A$ denote the $k$-submodule of $M$
consisting of the elements $m$ verifying $am = ma$ for all $a\in A$.

\remark{Remark 3.2.4} From Corollary~3.2.3, it follows immediately that if
$A$ is separable, then $\H^*(E,M) = \H^*\bigl(H,M^A\bigr)$ and if $A$ is
quasi-free, then there is a long exact sequence
$$
\align
\dots @>>> &\H^{n-2}(H,\H^1(A,M)) @>>>\H^n(H,\H^0(A,M)) @>>> \H^n(E,M)
@>>>\\
&\H^{n-1}(H,\H^1(A,M)) @>>> \H^{n+1}(H,\H^0(A,M)) @>>>\H^{n+1}(E,M)
@>>>\dots.
\endalign
$$
\endremark

\subhead 3.2.5. Separable subalgebras\endsubhead Let $S$ be a separable
subalgebra of $A$ and let $\wt{A}^r$ ($r\ge 0$) be as in 3.1.5. Suppose
$f(h,l)\in S$ for all $h,l\in H$. Using the fact that $f$ takes its values
in $S$ it is easy to see that $H$ acts on $\bigl(\Hom_{A^e}(A\ot_s\wt{A}^r
\ot_s A,M),b^*\bigr) = \bigl(\Hom_{S^e}( \wt{A}^r,M),b^*\bigr)$ via
$\bigl(\varphi\cdot h\bigr)(\wt{\ba}) = (1\# h^{(1)})^{-1}\varphi \bigl(
\wt{\ba}^{h^{(2)}}\bigr)(1\# h^{(3)})$.

\proclaim{Theorem 3.2.5.1} The Hochschild cohomology $\H^*(E,M)$, of $E$
with coefficients in $M$, is the cohomology of $H$ with coefficients in
$\bigl(\Hom_{S^e}(\wt{A}^r,M),b^*\bigr)$.

\endproclaim

\demo{Proof} It is similar to the proof of Theorem~3.1.5.1\qed
\enddemo

\subhead 3.2.6. An application to $\Ext_E^*$\endsubhead Let $k$ be a field,
$B$ an arbitrary $k$-algebra and $M$, $N$ two left $B$-modules. It is well
known that $\Ext_B^*(M,N)\simeq \H^*(B,\Hom_k(M,N))$ (here $\Hom_k(M,N)$ is
an $B$-bimodule via $(a\varphi b)(m) = a\varphi(bm)$). This fact and
Corollary~3.2.3 show that if $k$ is a field and $M$ and $N$ are left
$E$-modules, then there is a convergent spectral sequence
$$
E^{rs}_2 = \H^r(H,\Ext_A^s(M,N)) \Rightarrow \Ext_E^{r+s}(M,N).
$$
As a corollary we obtain that $\gl.dim(E) \le \gl.dim(A) + \gl.dim(H)$,
where $\gl.dim$ denotes the left global dimension. Note that this result
implies Maschke's Theorem for crossed product, as it was established in
\cite{B-M}.

\newpage

\head 4. The Cartan-Leray and Grothendieck spectral sequences \endhead
Assume that $E$ is a crossed product with invertible cocycle. In this case
another two spectral sequences converging to $\H_*(E,M)$ and with
$E^2$-term $\H_*(H,\H_*(A,M))$ can be considered. They are the Cartan-Leray
and the Grothendieck spectral sequences. The last one was introduced for
the more general setting of Galois extension in \cite{S}. In this Section
we recall these constructions and we prove that both coincide with the
Hochschild-Serre spectral sequence. Similar results are valid in the
cohomological setting.

\smallskip

Let $(\ov{H}^*\ot H,d_*)$ be the canonical resolution of $k$ as a right
$H$-module and $(Z_*,\partial_*) = (E\ot \ov{E}^*\ot E,b'_*)\ot
(\ov{H}^*\ot H,d_*)$. Consider $E\ot \ov{E}^r\ot E\ot \ov{H}^s\ot H$ as an
$E$-bimodule via
$$
(a\# l)(\bx\ot \bh)(b\# q) = ((a\# l)x_0\ot \bx_{1r}\ot x_{r+1}(b \#
q^{(1)}))\ot (\bh_{1s} \ot h_{s+1}q^{(2)}) ,
$$
where $\bx = x_0\ot\cdots\ot x_{r+1}$ and $\bh = h_1\ot\cdots\ot h_{s+1}$.
It is clear that
$$
E @<\mu << Z_0 @<\partial_1 << Z_1 @<\partial_2<< Z_2 @<\partial_3 << Z_3
@<\partial_4 << Z_4 @<\partial_5 << Z_5 @<\partial_6 << Z_6 @<\partial_7 <<
Z_7 @<\partial_8 << \dots,\tag3
$$
where $\mu((a_0\# h_0\ot a_1\# h_1)\ot l) = \ep(l)a_0a_1f(h_0^{(1)}
h_1^{(1)}) \# h_0^{(2)}h_1^{(2)}$, is a complex of $E$-bimodules. Moreover
$(3)$ is contractible as a complex of left $E$-modules, with contracting
homotopy $\zeta_n$ ($n\ge 0$) given by $\zeta_0(1_E) = 1_E\ot 1_E \ot 1_H$
and
$$
\zeta_{n+1}(\bby) = \cases - \bx\ot 1_E \ot \bh + (-1)^{n+1} x_0x_1\ot 1_E
\ot \bh\ot 1_H &\text{if $r = 0$} \\ (-1)^{r+1} \bx\ot 1_E \ot \bh
&\text{if $r>0$} \endcases,
$$
where $\bby = \bx \ot \bh$, with $\bx = x_0\ot\cdots \ot x_{r+1}$ and $\bh
= h_1\ot\cdots\ot h_{n-r+1}$. Since the map
$$
\tau\:E\ot \ov{E}^r\ot \ov{H}^s\ot H\ot E \to E\ot \ov{E}^r\ot E\ot
\ov{H}^s\ot H,
$$
given by $\tau(\bx_{0r}\ot \bh\ot x_{r+1}) = (\bx_{0r}\ot 1_E\ot
\bh)x_{r+1}$, is an isomorphism of $E$-bimodules (the inverse of $\tau$ is
the map $\bx_{0r}\ot a\# h \ot \bh\mapsto \bx_{0r} \ot \bh_{1s}\ot
h_{s+1}S^{-1}(h^{(2)})\ot a\# h^{(1)}$), $(Z_*,\partial_*)$ is a relative
projective resolution of $E$.

\bigskip

Let $M$ be an $E$-bimodule. The groups $M\ot_{E\ot A^{\op}} (E\ot
\ov{E}^r\ot E)$ are left $H$-modules via $h(m\ot \bx) = (1\# h^{(2)}) m \ot
\bx_{0r}\ot x_{r+1}(1\#h^{(1)})^{-1}$, where $\bx = x_0\ot\cdots\ot
x_{r+1}$. There is an isomorphism
$$
M \ot_{E^e} (Z_*,\partial_*) \simeq (\ov{H}^*\ot H,d_*)\ot_H (M\ot_{E\ot
A^{\op}} (E\ot \ov{E}^*\ot E,b'_*)).
$$
Let $F^i = \bigoplus_{j=0}^i (\ov{H}^j\ot H) \ot_H (M\ot_{A\ot E^{\op}}
E\ot\ov{E}^*\ot E)$. It is immediate that $F^0\sub F^1\sub F^2\sub F^3\sub
\dots$, is a filtration of the last complex. The spectral sequence
associate to this filtration converges to $\H_*(E,M)$ and has $E^2$-term
$\H_*(H,\H_*(A,M))$. This spectral sequence is called the homological
Cartan-Leray spectral sequence. Similarly the groups $\Hom_{E\ot A^{\op}}
(E\ot \ov{E}^r\ot E, M)$ are right $H$ modules via $f.h(\bx_{0,r+1}) =
f(\bx_{0r}\ot x_{r+1}(1\# h^{(1)})^{-1})(1\# h^{(2)})$ and there is an
isomorphism
$$
\Hom_{E^e}((Z_*,\partial_*),M) \simeq \Hom_H\bigl((\ov{H}^*\ot
H,d_*),\Hom_{E\ot A^{\op}}((E\ot \ov{E}^*\ot E,b'_*),M)\bigr).
$$
This complex has a filtration $F_0\sup F_1\sup F_2\sup F_3\sup F_4\sup
\dots$, defined by $F_i^n = \bigoplus_{j\ge i} \Hom_H\bigl(\ov{H}^j\ot
H,\Hom_{E\ot A^{\op}}(E\ot \ov{E}^{n-j}\ot E,M)\bigr)$. The spectral
sequence associate to this filtration converges to $\H^*(E,M)$ and has
$E^2$-term $\H^*(H,\H^*(A,M))$. This spectral sequence is called the
cohomological Cartan-Leray spectral sequence.

\bigskip

Let $\Phi_*\: (E\ot \ov{E}^*\ot E, b'_*)\to (Z_*,\partial_*)$ and $\Psi_*\:
(Z_*,\partial_*) \to (E\ot \ov{E}^*\ot E, b'_*)$ be the morphisms of
$E$-bimodule complexes, recursively defined by
$$
\align
&\Phi_0(x\ot 1_E) = x\ot 1_E\ot 1_H, \quad \Psi_0(x\ot 1_E\ot h) =
\ep(h)x\ot 1_E,\\
& \Phi_{n+1}(\bx\ot 1_E) = \zeta_{n+1} \circ \Phi_n\circ b'_{n+1}(\bx\ot
1_E)\,\text{ for $\bx\in E\ot\ov{E}^{n+1}$,}\\
& \Psi_{n+1}(\bx\ot 1_E\ot \bh) = \xi_{n+1} \circ \psi_n \circ
\partial_{n+1} (\bx\ot 1_E\ot \bh)\,\text{ for $\bx\in E\ot\ov{E}^r$,
$\bh\in \ov{H}^{n+1-r}\!\ot\! H$.}
\endalign
$$

\proclaim{Proposition 4.1} It is hold that $\Psi_*\circ\Phi_* = id_*$ and
that $\Phi_*\circ\Psi_*$ is homotopically equivalent to the identity map.
The homotopy $\Phi_*\circ \Psi_*  @>\Om_{*+1}>> id_*$ is recursively
defined by $\Om_1(x\ot 1_E\ot h) = x\ot 1_E\ot h\ot 1_H$ and
$$
\Om_{n+1}(\bx\ot 1_E\ot \bh) = \zeta_{n+1}\circ (\Phi_n \circ \Psi_n - id -
\Om_n \circ \partial_n)(\bx\ot 1_E\ot \bh),
$$
for $\bx = x_0\ot\cdots\ot x_r$ and $\bh = h_1\ot \cdots \ot h_{n+1-r}$.

\endproclaim

\demo{Proof} It is easy to see that $\Phi_*$ and $\Psi_*$ are morphisms of
complexes. Arguing as in Proposition 1.2.1 we get that $\Om_{*+1}$ is an
homotopy from $\Phi_* \circ \Psi_*$ to the identity map. It remains to
prove that $\Psi_*\circ \Phi_* = id_*$. It is clear that $\Psi_0\circ\Phi_0
= id_0$. Assume that $\Psi_n\circ\Phi_n = id_n$. Since $\Phi_{n+1}(E\ot
\ov{E}^n\ot k)\sub \sum_{r=0}^{n+1} E\ot \ov{E}^r\ot k\ot \ov{H}^{n+1-r}\ot
H$, we have that on $E\ot \ov{E}^n\ot k$
$$
\align
\Psi_{n+1}\circ \Phi_{n+1}& = \xi_{n+1}\circ\Psi_n \circ
\partial_{n+1}\circ \Phi_{n+1} = \xi_{n+1}\circ\Psi_n\circ
\partial_{n+1}\circ \zeta_{n+1}\circ\Phi_n\circ b'_{n+1}\\
& = \xi_{n+1}\circ\Psi_n\circ\Phi_n\circ b'_{n+1} - \xi_{n+1}\circ\Psi_n
\circ \zeta_n\circ \partial_n\circ\Phi_n \circ b'_{n+1} = \xi_{n+1}\circ
b'_{n+1} = id_{n+1}\qed
\endalign
$$

\enddemo

Next, we consider the normalized Hochschild resolution $(E\ot \ov{E}^*\ot
E,b'_*)$ filtered as in Proposition~1.2.2 and the resolution $(Z_*,
\partial_*)$ filtered by $F^0_* \sub F^1_* \sub F^2_* \sub \dots$, where
$F^i_* = \bigoplus_{j=0}^i (E\ot \ov{E}^{n-j}\ot E)\ot (\ov{H}^j\ot H)$.

\proclaim{Proposition 4.2} We have that
$$
\align
\Phi_n(a_0\# h_0&\ot\cdots \ot a_{n+1}\# h_{n+1}) = \sum_{j=0}^n
(-1)^{j(n+1)}(a_0\# h_0)(a_1\# h_1^{(1)}) \dots (a_j\# h_j^{(1)})\\
&\ot (a_{j+1}\# h_{j+1}^{(1)}) \ot \cdots \ot (a_{n+1}\# h_{n+1}^{(1)})\!
\ot h_1^{(2)}\ot\cdots \ot h_j^{(2)}\!\ot h_{j+1}^{(2)}\cdots h_{n+1}^{(2)}.
\endalign
$$
Consequently the map $\Phi_*$ preserve filtrations.
\endproclaim

\demo{Proof} It follows by induction on $n$, using the recursive definition
of $\Phi_*$\qed
\enddemo

\proclaim{Proposition 4.3} The map $\Phi_*$ induces an homotopy equivalence
of $E$-bimodule complexes between the graded complexes associated to the
filtrations of $(\B_*(E),b'_*)$ and $(\B_*(E),b'_*)\ot (\ov{H}^* \ot
H,d_*)$.
\endproclaim

\demo{Proof} Note that
$$
\align
& \frac{F^s(X_*,d_*)}{F^{s-1}(X_*,d_*)} = (X_{*s},d^0_{*s}) = (E\ot
\ov{H}^s\ot \ov{A}^*\ot E,d^0_{*s}), \\
&\frac{F^s(Z_*,\partial_*)}{F^{s-1}(Z_*,\partial_*)} = (\B_*(E),b'_*)\ot
\ov{H}^s\ot H,
\endalign
$$
where $d^0_{*,s}$ is the boundary map introduced in Subsection~1.1. By
Proposition~1.2.2 it suffices to check that $\varPhi_* = \Phi_*\circ
\phi_*$ induces an homotopy equivalence $\wt{\varPhi}^s_*$ of
$E$-bimodules complexes, from $(E\ot \ov{H}^s\ot \ov{A}^*\ot E,d^0_{*,s})$
to $(\B_*(E),b'_*)\ot \ov{H}^s\ot H$. Let $Y_s$ and $\mu_s$ be as in
Subsection~1.1 and $\wt{Y}_s = E\ot \ov{H}^s \ot H$ endowed with the structure
of $E$-bimodule given by $x_0(x_1\ot \bh)x_2 = x_0x_1x_2 \ot\bh$, where
$\bh = h_0\ot\cdots\ot h_{s+1}$. Consider the diagram
$$
\CD
Y_s @<\mu_s << E\ot \ov{H}^s \ot E @<d^0_{1s}<< E\ot \ov{H}^s\ot \ov{A}\ot
E @<d^0_{2s}<< \dots\\
@VV\wt{\varPhi}^s V @VV\wt{\varPhi}^s_0 V @VV\wt{\varPhi}^s_1V \\
\wt{Y}_s @<\wt{\mu}_s << E\ot E\ot \ov{H}^s\ot H @<b'_1<< E\ot\ov{E}\ot E
\ot \ov{H}^s\ot H @<b'_2<< \dots,
\endCD \tag4
$$
where $\wt{\mu}_s((x_0\ot x_1)\ot \bh) = x_0x_1\ot \bh$ and
$\wt{\varPhi}^s(x \ot \bh) = x(1\# h_1^{(1)}) \cdots (1\# h_{s+1}^{(1)})\ot
\bh^{(2)}$. We assert that $\wt{\varPhi}^s_0 (\bx) = 1_E\ot (1\#
h_1^{(1)})\cdots (1\# h_{s}^{(1)})\ot \bh^{(2)}\ot 1_H$, where $\bx =
1_E\ot \bh\ot 1_E$, with $\bh = h_1\ot\cdots\ot h_s$. To prove this it
suffices to check that
$$
\Phi_s\circ \phi_s (\bx) \in 1_E\ot (1\# h_1^{(1)})\cdots (1\# h_s^{(1)})
\ot \bh^{(2)}\ot 1_H + F_{s-1},
$$
which follows by induction on $s$, using that $\Phi_s\circ \phi_s (\bx) =
\zeta_s \circ \Phi_{s-1}\circ \phi_{s-1} \circ d_s(\bx)$. Now, it is
immediate that $\wt{\mu}_s \circ \wt{\varPhi}^s_0 = \wt{\varPhi}^s \circ
\mu_s$. Since $\wt{\varPhi}^s$ is an isomorphism and the rows of $(4)$ are
relative projective resolutions of $Y_s$ and $\wt{Y}_s$ respectively, it
follows that $\wt{\varPhi}^s_*$ is an homotopy equivalence\qed

\enddemo

\proclaim{Corollary 4.4} The (co)homological Cartan-Leray spectral sequence
is isomorphic to the (co)homological Hochschild-Serre spectral sequence.
\endproclaim

\subheading{4.5 The Grothendieck spectral sequence} If $M$ is an
$E$-bimodule, then the group $H_0(A,M) = M/[A,M]$ is a left $H$-module via
$h\cdot \ov{m} = \ov{(1\# h^{(2)}) m (1\# h^{(1)})^{-1}}$, where the
$\ov{m}$ denotes the class of $m$ in $M/[A,M]$. Let us consider the
functors $M\mapsto H_0(E,M)$ from the category of $E$-bimodules to the
category of $k$-modules, $M\mapsto H_0(A,M)$ from the category of
$E$-bimodules to the category of left $H$-modules and $M\mapsto H_0(H,M)$
from the category of left $H$-modules to the category of $k$-modules. It is
easy to see that $H_0(E,M) = H_0(H,H_0(A,M))$ and that if $M$ is a
relatively projective $E^e/E^{\op}$-module, then $H_0(A,M)$ is a relatively
projective $H/k$-module. In fact, if $M = E\ot N$, then the map $h\ot n
\mapsto \ov{1\# h^{(2)}\ot n(1\# h^{(1)})^{-1}}$ is an isomorphism of left
$H$-modules from $H\ot N$ to $H_0(A,M)$. Thus we have a Grothendieck
spectral sequence
$$
E^2_{rs} = \H_s(H,\H_r(A,M)) \to \H_{r+s}(E,M)).
$$
We assert that the Grothendieck spectral sequence and the Cartan-Leray
spectral sequence coincide. To prove this we use a concrete construction of
the Grothendieck spectral sequence. Let $(P_*,\partial_*) = (M\ot
\ov{E}^*\ot E,b'_*)$ be the normalized canonical resolution of $M$ as a
right $E^e$-module. Let us write $(\ov{P}_*,\ov{\partial}_*) =
(P_*,\partial_*)\ot_{A^e} A$. Consider the double complex
$$
C_{**}:=
\CD
\vdots &&\vdots && \vdots \\
@VVV   @VVV  @VVV \\
H\ot_H \ov{P}_1 @<<< \ov{H}\ot H\ot_H \ov{P}_1 @<<< \ov{H}^2\ot
H\ot_H\ov{P}_1 @<<< \dots \\
@VVV   @VVV  @VVV \\
H\ot_H\ov{P}_0 @<<< \ov{H}\ot H\ot_H\ov{P}_0 @<<< \ov{H}^2\ot
H\ot_H\ov{P}_0 @<<< \dots,
\endCD
$$
whose $r$-th column is $(-1)^r$ times $\ov{H}^r\ot H\ot_H (\ov{P}_*,
\ov{\partial}_*)$ and whose $s$-th row is the canonical complex $(\ov{H}^*
\ot H\ot_H \ov{P}_s , d_*)$ giving the homology $H_*(H,\ov{P}_s)$ of $k$ as
a trivial right $H$-module with coefficients in $\ov{P}_s$. By definition,
the Grothendieck spectral sequence is the spectral sequence associate to
the filtrations by columns of $C_{**}$. Since $C_{**} \simeq (\ov{H}^*\ot
H,d_*)\ot_H (M\ot_{E\ot A^{\op}} (E\ot \ov{E}^*\ot E,b'_*))$ as filtered
complexes, the homological Cartan-Leray and the Grothendieck spectral
sequence coincide. The same is valid in the cohomological setting.

\head Appendix A \endhead
Let $R\to S$ be an unitary ring map and let $N$ be a left $S$-module. In
this section, under suitable conditions, we construct a projective relative
resolution of $N$. We need this result (with $R = E$, $S = E^e$ and $N =
E$) to complete the proof of Theorem~1.1.1. The general case considered here
simplifies the notation and enables us to consider other cases, for
instance algebras of groups having particular resolutions.

\smallskip

Let us consider a diagram of left $S$-modules and $S$-module maps
$$
\CD
\vdots\\
@VV\partial_2V \\
Y_1 @<\mu_1<<  X_{01} @<d^0_{11}<<  X_{11} @<d^0_{21}<< \dots \\
@VV\partial_1V \\
Y_0 @<\mu_0<<  X_{00} @<d^0_{10}<<  X_{10} @<d^0_{20}<< \dots,
\endCD
$$
such that:

\smallskip

\item{a)} The column and the rows are chain complexes.

\smallskip

\item{b)} For each $r,s\ge 0$ we have a left $R$-module $\ov{X}_{rs}$ and
$S$-module maps
$$
s_{rs}\:X_{rs}\to S\ot \ov{X}_{rs}\qquad\text{and}\qquad \pi_{rs}\:S\ot
\ov{X}_{rs}\to X_{rs}
$$
verifying $\pi_{rs}\circ s_{rs} = id$.

\smallskip

\item{c)} Each row is contractible as a complex of $R$-modules, with a
chain contracting homotopy $\si^0_{0s}\: Y_s \to X_{0s}$ and $\si^0_{r+1,s}
\: X_{rs}\to X_{r+1,s}$ ($r\ge 0$).

\medskip

\ni We are going to modify this diagram by adding $S$-module maps
$$
d^l_{rs}\:X_{rs} @>>> X_{r+l-1,s-l}\qquad\text{($r,s\ge 0$ and $1\le l\le
s$).}
$$
Let $X_n = \bigoplus_{r+s=n} X_{rs}$ and $d_n = \sum_{r+s=n} \sum_{l=0\atop
r+l>0}^s d^l_{rs}$ ($n\ge 1$). Consider the maps $\mu'_n\: X_n \to Y_n$
($n\ge 0$), given by:
$$
\mu'_n (x) = \cases \mu_n(x) &\text{for $x\in X_{0n}$}\\ 0 &\text{for
$x\in X_{r,n-r}$ with $r>0$.} \endcases
$$
We define the arrows $d^l_{rs}$ in such a way that $(X_*,d_*)$ becomes a
chain complex of $S$-modules and $\mu'_*\: (X_*,d_*) \to (Y_*,-\partial_*)$
becomes a chain homotopy equivalence of complexes of $R$-modules. In fact,
we are going to build $R$-module morphisms
$$
\si^l_{l,s-l} \: Y_s \to X_{l,s-l}\quad\!\!\text{and}\quad\!\!
\si^l_{r+l+1,s-l} \: X_{rs} \to X_{r+l+1,s-l}\quad\!\!\text{($r,s\ge 0$ and
$1\le l\le s$),}
$$
satisfying the following:

\proclaim{Theorem A.1} Let $\C_*(\mu'_*)$ be the mapping cone of $\mu'_*$,
that is, $\C_*(\mu'_*) = (C_*,\de_*)$, where $C_n = Y_n\oplus X_{n-1}$ and
$\de_n(y_n,x_{n-1}) = \bigl(-\partial(y_n) - \mu'_{n-1}(x_{n-1}),-d_{n-1}
(x_{n-1})\bigr)$. The family of $R$-module maps $\si_{n+1}\: \C_n(\mu'_*)
\to \C_{n+1}(\mu'_*)$ ($n\ge 0$), defined by:
$$
\si_{n+1} = -\sum_{r+s=n-1\atop r\ge -1} \sum_{l=0}^s \si^l_{r+l+1,s-l},
$$
is a chain contracting homotopy of $\C_*(\mu'_*)$.
\endproclaim

\proclaim{Corollary A.2} Let $N$ be a left $S$-module. If there is a
$S$-module map $\wt{\mu}\: Y_0 \to N$, such that
$$
N @<\wt{\mu}<< Y_0  @<\partial_1<< Y_1  @<\partial_2 << Y_2 @<\partial_3<<
Y_3 @<\partial_4<< Y_4 @<\partial_5<< Y_6 @<\partial_7<<\dots \tag *
$$
is contractible as a complex of left $R$-modules, then
$$
N @<\mu<< X_0  @<d_1<< X_1  @<d_2<< X_2  @<d_3<< X_3 @<d_4<< X_4 @<d_5<<
X_5 @<d_6<< X_6 @<d_7<<\dots,  \tag **
$$
where $\mu = \wt{\mu}\circ \mu_0$, is a relative projective resolution.
Moreover, if $\si_0^{-1}\: N \to Y_0$, $\si_{n+1}^{-1}\:Y_n \to Y_{n+1}$
($n\ge 0$) is a chain contracting homotopy of (*), then we obtain a chain
contracting homotopy $\ov{\si}_0\: N \to X_0$, $\ov{\si}_{n+1}\:X_n \to
X_{n+1}$ ($n\ge 0$) of (**), defining $\ov{\si}_0 = \si_{00}^0 \circ
\si_0^{-1}$ and
$$
\ov{\si}_{n+1} = - \sum_{l=0}^{n+1} \si_{l,n-l+1}^l \circ
\si_{n+1}^{-1}\circ \mu_n + \sum_{r+s=n} \sum_{l=0}^s \si_{r+l+1,s-l}^l.
$$
\endproclaim

\demo{Proof} Write
$$
\wt{\si}_n = \sum_{r+s=n-1\atop r\ge 0} \sum_{l=0}^s \si_{r+l+1,s-l}^l\quad
\text{($n\ge 1$)}\quad\text{and}\quad\wh{\si}_n = \sum_{l=0}^n\si_{l,n-l}^l
\quad\text{($n\ge 0$).}
$$
From Theorem A.1, we have
$$
\wh{\si}_n \circ \partial_{n+1} = \sum_{l=0}^n\si^l_{l,n-l}\circ
\partial_{n+1} = - \sum_{l=0}^n \sum_{i=0}^{l+1} d^{l+1-i}_{i,n+1-i}\circ
\si^i_{i,n+1-i} = - d_{n+1}\circ \wh{\si}_{n+1}. \tag *n
$$
It is clear that $\mu \circ \ov{\si}_0 = id$. Moreover
$$
\align
\ov{\si}_0 \circ \mu & = \si_{00}^0 \circ \si_0^{-1} \circ \wt{\mu}\circ
\mu_0 = \si_{00}^0 \circ \mu_0 - \si_{00}^0 \circ \partial_1 \circ
\si_1^{-1}\circ \mu_0 \\
& = id - d_{10}^0\circ\si_{10}^0 + d_{01}^1\circ\si_{01}^0\circ\si_1^{-1}
\circ \mu_0 + d_{10}^0\circ \si_{10}^1\circ \si_1^{-1}\circ \mu_0,
\endalign
$$
where the last equality follows from (*0). Now, let $n\ge 1$. Take $x\in
X_{r,n-r}$. If $r\ge 1$, then the equality $(0,x) = \de_{n+2}\circ
\si_{n+2}(0,x) + \si_{n+1}\circ \de_{n+1}(0,x)$ implies that $x =
d_{n+1}\circ \wt{\si}_{n+1}(x) + \wt{\si}_n\circ d_n (x)$. Hence, we can
suppose $r = 0$. Then, from $(0,x) = \de_{n+2}\circ \si_{n+2}(0,x) +
\si_{n+1}\circ \de_{n+1}(0,x)$, we get
$$
\align
x & = d_{n+1}\!\!\circ \wt{\si}_{n+1}(x) + \wt{\si}_n\!\!\circ d_n(x) +
\wh{\si}_n \!\!\circ \mu_n(x) \\
& = d_{n+1}\!\!\circ \wt{\si}_{n+1}(x) + \wt{\si}_n\!\!\circ d_n(x) +
\wh{\si}_n \!\!\circ \si_n^{-1}\!\!\circ \partial_n\!\!\circ \mu_n(x) +
\wh{\si}_n \!\!\circ \partial_{n+1}\!\!\circ \si_{n+1}^{-1}\!\!\circ
\mu_n(x) \\
& = d_{n+1}\!\!\circ \wt{\si}_{n+1}(x) + \wt{\si}_n\!\!\circ d_n(x) -
\wh{\si}_n \!\!\circ \si_n^{-1}\!\!\circ \mu_{n-1}\!\!\circ d_n(x) +
\wh{\si}_n \!\!\circ \partial_{n+1}\!\!\circ \si_{n+1}^{-1}\!\!\circ
\mu_n(x)\\
& = d_{n+1}\!\!\circ \wt{\si}_{n+1}(x) + \wt{\si}_n\!\!\circ d_n(x) -
\wh{\si}_n \!\!\circ \si_n^{-1}\!\!\circ \mu_{n-1}\!\!\circ d_n(x) -
d_{n+1}\!\!\circ\!\wh{\si}_{n+1} \!\!\circ  \si_{n+1}^{-1}\!\!\circ
\mu_n(x),
\endalign
$$
where the last equality follows from (*n)\qed
\enddemo

\medskip

Next we define the morphisms $d^l_{rs}$ and we prove that $(X_*,d_*)$ is a
chain complex.

\proclaim{Definition A.3} We define the $S$-module maps $d^l_{rs}\: X_{rs}
\to X_{r+l-1,s-l}$ ($r\ge 0$ and $1\le l\le s$), recursively by $d^l_{rs} =
\ov{d}^l_{rs}\circ s_{rs}$, where $\ov{d}^l_{rs}\: S\ot \ov{X}_{rs}\to
X_{r+l-1,s-l}$ ($r\ge 0$ and $1\le l\le s$) is the $S$-module map defined
by
$$
\ov{d}^l_{rs}(\bx) = \cases
- \si_{0,s-1}^0\circ \partial_s \circ \mu_s \circ \pi_{0s} (\bx)
&\text{if $r=0$ and $l=1$,}\\
- \sum_{j=1}^{l-1} \si_{l-1,s-l}^0\circ d^{l-j}_{j-1,s-j} \circ
d^j_{0s}\circ \pi_{0s} (\bx)&\text{if $r=0$ and $1<l\le s$,}\\
- \sum_{j=0}^{l-1} \si_{r+l-1,s-l}^0\circ d^{l-j}_{r+j-1,s-j} \circ
d^j_{rs}\circ \pi_{rs} (\bx) &\text{if $r>0$,}
\endcases
$$
for each $\bx = 1\ot \ov{\bx}\in S\ot\ov{X}_{rs}$.
\endproclaim

\proclaim{Proposition A.4} We have $\mu_{s-1}\circ d^1_{0s} = -
\partial_s\circ \mu_s$ and
$$
d^0_{r+l-1,s-l} \circ d^l_{rs} = \cases - \sum_{j=1}^{l-1}
d^{l-j}_{j-1,s-j} \circ d^j_{0s} &\text{if $r=0$ and $1<l\le s$}\\
- \sum_{j=0}^{l-1} d^{l-j}_{r+j-1,s-j} \circ d^j_{rs}&\text{if $r>0$ and
$1\le l\le s$.}
\endcases
$$
Consequently $(X_*,d_*)$ is a chain complex.

\endproclaim

\demo{Proof} We prove the proposition by induction on $l$ and $r$. To
simplify the expressions we put $d^0_{0s}:=\mu_s$, $d^1_{-1,s}:=\partial_s$
and $d^l_{-1,s}:= 0$ for all $l>1$. Moreover to abbreviate we do not write
the subindices. Let $\bx = 1\ot \ov{\bx}$ with $\ov{\bx}\in \ov{X}_{0s}$.
Since $\ov{d}^1_0(\bx) = - \si^0\circ d^1\circ d^0\circ \pi(\bx)$, we have
$d^0\circ \ov{d}^1(\bx) = - d^0 \circ \si^0\circ d^1\circ d^0\circ \pi(\bx)
= - d^1\circ d^0\circ \pi(\bx)$, which implies $d^0\circ d^1 = - d^1\circ
d^0$. Let $l+r>1$ and suppose the result is valid for $d^j_{p*}$ with $j<l$
or $j=l$ and $p<r$. Let $\bx = 1\ot\ov{\bx}$ with $\ov{\bx}\in
\ov{X}_{rs}$. Since $ \ov{d}^l(\bx) = - \sum_{j=0}^{l-1} \si^0\circ d^{l-j}
\circ d^j\circ \pi(\bx)$, then
$$
d^0\circ \ov{d}^l(\bx) = - \sum_{j=0}^{l-1} d^0\circ \si^0 \circ d^{l-j}
\circ d^j\circ \pi(\bx) = - \sum_{j=0}^{l-1} d^{l-j}\circ d^j\circ \pi(\bx)
+ \sum_{j=0}^{l-1} \si^0\circ d^0 \circ d^{l-j} \circ d^j\circ \pi(\bx).
$$
Applying first the inductive hypothesis to $d^0\circ d^{l-j}$ with ($0\le
j<l$) and then to $d^0\circ d^j$ with ($0<j< l$), we obtain:
$$
\align
d^0\circ \ov{d}^l(\bx) &= - \sum_{j=0}^{l-1} d^{l-j} \circ d^j\circ
\pi(\bx) - \sum_{j=0}^{l-1} \sum_{i=0}^{l-j-1} \si^0 \circ d^{l-j-i} \circ
d^i \circ d^j\circ \pi (\bx)\\
&= - \sum_{j=0}^{l-1} d^{l-j} \circ d^j\circ \pi(\bx) - \sum_{j=0}^{l-2}
\sum_{i=1}^{l-j-1} \si^0\circ d^{l-j-i} \circ d^i \circ d^j\circ \pi
(\bx)\\
& + \sum_{j=1}^{l-1} \sum_{h=0}^{j-1} \si^0\circ d^{l-j} \circ d^{j-h}\circ
d^h \circ \pi(\bx) = - \sum_{j=0}^{l-1} d^{l-j} \circ d^j\circ \pi(\bx).
\endalign
$$
The desired equality follows immediately from this fact\qed
\enddemo

It is immediate that $\mu'_*\: (X_*,d_*) \to (Y_*,-\partial_*)$ is a
morphism of $S$-module chain complexes. Next, we construct the chain
contracting homotopy of $\C_*(\mu'_*)$.

\proclaim{Definition A.5}   We define $\si^l_{l,s-l}\: Y_s \to X_{l,s-l}$
and $\si^l_{r+l+1,s-l}\: X_{rs} \to X_{r+l+1,s-l}$ ($0<l\le s$, $r\ge 0$),
recursively by:
$$
\si^l_{r+l+1,s-l} = - \sum_{i=0}^{l-1} \si^0_{r+l+1,s-l} \circ d^{l-i}_{
r+i+1,s-i} \circ\si^i_{r+i+1,s-i}\quad\text{($0<l\le s$ and $r\ge -1$)}.
$$
\endproclaim

\subhead Proof of Theorem~A.1 \endsubhead To simplify the expressions we
put $d^0_{-1,s}:= 0$, $d^0_{0s}:=\mu_s$, $d^1_{-1,s}:=\partial_s$ and
$d^l_{-1,s}:= 0$ for all $l>1$. Because of the definitions of $d_*$ and
$\si_*$, it suffices to check that $\si^0_{rs}\circ d^0_{rs} + d^0_{r+1,s}
\circ \si^0_{r+1,s} = id$ and
$$
\sum_{i=0}^l \si^{l-i}_{r+l,s-l}\circ d^i_{rs} + \sum_{i=0}^l
d^{l-i}_{r+i+1,s-i}\circ \si^i_{r+i+1,s-i} = 0\quad\text{for $l>0$,}
$$
where we put $d^0_{-1,s} = 0$. The first formula simply says that $\si^0_*$
is a chain contracting homotopy of $d^0_*$. Let us see the second one. To
abbreviate we do not write the subindices. From the definition of $\si^l$ we
have:
$$
d^0\circ \si^l = - \sum_{i=0}^{l-1} d^0\circ \si^0\circ d^{l-i}\circ \si^i
= \sum_{i=0}^{l-1} \si^0\circ d^0\circ d^{l-i}\circ \si^i - \sum_{i=0}^{l-1}
d^{l-i}\circ \si^i.
$$
Consequently
$$
\sum_{i=0}^l \si^{l-i}\circ d^i + \sum_{i=0}^l d^{l-i}\circ \si^i =
\sum_{i=0}^l \si^{l-i} \circ d^i + \sum_{i=0}^{l-1} \si^0\circ d^0\circ
d^{l-i}\circ \si^i.
$$
Then, it suffices to prove that the term appearing on the right side of the
equality is zero. We prove this by induction on $l$. For $l=1$ we have:
$$
\si^0\circ d^0\circ d^1\circ \si^0 = - \si^0\circ d^1\circ d^0\circ \si^0 =
\si^0\circ d^1\circ \si^0\circ d^0 - \si^0\circ d^1  =  - \si^1\circ d^0 -
\si^0\circ d^1.
$$
Suppose $l>1$. From Proposition~A.5,
$$
\sum_{i=0}^{l-1} \si^0\circ d^0\circ d^{l-i}\circ \si^i = -
\sum_{i=0}^{l-1} \sum_{j=0}^{l-i-1} \si^0\circ d^{l-i-j}\circ d^j\circ
\si^i = - \sum_{h=0}^{l-1} \sum_{i=0}^h \si^0\circ d^{l-h}\circ
d^{h-i}\circ \si^i.
$$
So, applying the inductive hypothesis to $\sum_{i=0}^h d^{h-i}\circ \si^i$
($h\ge 0$), we obtain
$$
\align
\sum_{i=0}^{l-1} \si^0\circ d^0\circ d^{l-i}\circ \si^i &= \sum_{h=0}^{l-1}
\sum_{i=0}^h \si^0\circ d^{l-h}\circ \si^{h-i}\circ d^i - \si^0 \circ d^l\\
&= \sum_{i=0}^{l-1} \sum_{j=0}^{l-i-1} \si^0\circ d^{l-i-j}\circ \si^j
\circ d^i - \si^0\circ d^l\\
& = - \sum_{i=0}^l \si^{l-i}\circ d^i\qed
\endalign
$$

\head Appendix B\endhead
In this appendix we compute explicitly the maps $d_{rs}^l$ introduced in
Section~1, completing the results of Theorem~1.1.3.

\definition{Definition B.1} Given $\bh= h_1\ot\cdots\ot h_l\in \ov{H}^l$,
we define $F_0^{(l)}(\bh)$, recursively by:
$$
\align
& F_0^{(2)}(\bh) = -f(h_1,h_2), \\
& F_0^{(l+1)}(\bh) = \sum_{j=1}^l (-1)^j f(h_j^{(1)},h_{j+1}
^{(1)})^{\ov{\bh_{1,j-1}^{(1)}}}\ot F^{(l)}(\bh^{j(2)}),
\endalign
$$
where $\bh^{j(2)} = \bh_{1,j-1}^{(2)}\ot \h_j^{(2)} h_{j+1}^{(2)} \ot
\bh_{j+2,l+1}$. For instance, we have
$$
\align
F_0^{(3)}(\bh) = & f(h_1^{(1)}\!,h_2^{(1)})\ot f(\fh_{12}^{(2)}\!,h_3) -
f(h_2^{(1)}\!,h_3^{(1)})^{h_1^{(1)}} \ot f(h_1^{(2)}\!,\fh_{23}^{(2)})\\
\intertext{and}
F_0^{(4)}(\bh) = &  - f(h_1^{(1)}\!,h_2^{(1)})\ot f(\fh_{12}^{(2)}\!,
h_3^{(1)}) \ot f(\fh_{12}^{(3)}h_3^{(2)}\!,h_4)\\
& + f(h_1^{(1)}\!,h_2^{(1)})\ot f(h_3^{(1)}\!,h_4^{(1)})^{\fh_{12}^{(2)}}
\ot f(\fh_{12}^{(2)}\!,\fh_{34}^{(2)}) \\
& + f(h_2^{(1)}\!,h_3^{(1)})^{h_1^{(1)}}\ot f(h_1^{(2)}\!,\fh_{23}^{(2)})
\ot f(\fh_{13}^{(3)}\!,h_4) \\
& - f(h_2^{(1)}\!,h_3^{(1)})^{h_1^{(1)}}\ot f(\fh_{23}^{(2)}\!,h_4^{(1)})
^{h_1^{(2)}} \ot f(h_1^{(3)}\!,\fh_{23}^{(3)}h_4^{(2)}) \\
& - f(h_3^{(1)}\!,h_4^{(1)})^{\ov{\bh_{12}^{(1)}}}\ot f(h_1^{(2)}\!,
h_2^{(2)}) \ot f(\fh_{12}^{(3)}\!,\fh_{34}^{(2)})\\
& + f(h_3^{(1)}\!,h_4^{(1)})^{\ov{\bh_{12}^{(1)}}}\ot f(h_2^{(2)}\!,
\fh_{34}^{(2)})^{h_1^{(2)}} \ot f(h_1^{(3)}\!,\fh_{24}^{(3)}).
\endalign
$$
\enddefinition

\medskip

For the following definition we adopt the convention that $\ba_{10} =
\ba_{r+1,r} = 1_k\in k$.

\definition{Definition B.2} Given $\bh= h_1\ot\cdots\ot h_l\in \ov{H}^l$
and $\ba = a_1\ot\cdots\ot a_r\in \ov{A}^r$, we define $F_r^{(l)} (\bh\ot
\ba)$, recursively by:
$$
\align
& F_r^{(2)}(\bh\ot \ba) = \sum_{i=0}^r (-1)^{i+1}
\ba_{1i}^{\ov{\bh_{12}^{(1)} }}\ot f(h_1^{(2)},h_2^{(2)})\ot
\ba_{i+1,r}^{\fh_{12}^{(3)}},\\
& F_r^{(l+1)}(\bh\ot \ba ) = \sum_{j=1}^l\sum_{i=0}^r (-1)^{il+j}
\ba_{1i}^{\ov{\bh_{1,l+1}^{(1)}}}\ot f(h_j^{(2)},h_{j+1}^{(2)})^{
\ov{\bh_{1,j-1}^{(2)}}} \ot F_{r-i}^{(l)}(\bh^{j(3)} \ot \ba_{i+1,r}),
\endalign
$$
where  $\bh^{j(3)} = \bh_{1,j-1}^{(3)}\ot \h_j^{(3)} h_{j+1}^{(3)} \ot
\bh_{j+2,l+1}^{(2)}$ and $F_0^{(l)}( \bh^{j(3)}\ot \ba_{r+1,r}) =
F_0^{(l)}(\bh^{j(3)})$. For instance, we have
$$
\align
& F_r^{(3)}(\bh\ot\ba)= \sum_{0\le i\le j\le r}(-1)^{i+j} \ba_{1i}^{
\ov{\bh_{13}^{(1)}}}\ot f(h_1^{(2)}\!,h_2^{(2)})\ot \ba_{i+1,j}^{\ov{
\bh_{13}^{1(3)}}}\ot f(\fh_{12}^{(4)}\!,h_3^{(3)})\ot \ba_{j+1,r}^{\fh_{12}
^{(5)}h_3^{(4)}}\\
&\phantom{F_r^{(3)}(\bh} + \sum_{0\le i\le j\le r} (-1)^{i+j+1}
\ba_{1i}^{\ov{\bh_{13}^{(1)}}}\ot f(h_2^{(2)}\!,h_3^{(2)})^{h_1^{(2)}}\ot
\ba_{i+1,j}^{\ov{\bh_{13}^{2(3)}}} \ot f(h_1^{(4)}\!,\fh_{23}^{(4)})\ot
\ba_{j+1,r}^{\fh_{13}^{(5)}}.
\endalign
$$
\enddefinition

\smallskip

We set $F_0^{(1)}(h_s) = 1_k\in k$, $F_r^{(1)}(h_s\ot \ba) = \ba^{h_s}$ and
$F_0^{(l)} (\bh_{s-l-1,s}\ot 1_k) = F_0^{(l)}(\bh_{s-l-1,s})$. Moreover, to
abbreviate we write $F^{(l)}(\bh) = F_0^{(l)}(\bh)$ and
$F^{(l)}\!\!\left(\!\smallmatrix \ba\\ \bh\endsmallmatrix \! \right) =
F_r^{(l)}( \bh\ot \ba)$.

\proclaim{Lemma B.3} Let $\ba = a_1\ot\cdots\ot a_r$ and $\bh_{s-l,s} =
h_{s-l}\ot\cdots\ot h_l$. We have:
$$
\align
F^{(l+1)}(\bh_{s-l,s}) &= \sum_{i=1}^l (-1)^i F^{(l-i+1)}\!\!\left(\!
\smallmatrix F^{(i)}(\bh_{s-i+1,s}^{(1)})\\ \vspace{1pt} \bh_{s-l,s-i}
^{(1)} \endsmallmatrix \!\!\!\right) \ot f(\fh_{s-l,s-i}^{(2)},
\fh_{s-i+1,s}^{(2)})\\
\intertext{and}
F^{(l+1)}\!\!\left(\!\smallmatrix \ba\\ \bh_{s-l,s}
\endsmallmatrix\!\!\right) & = F^{(l+1)}\!\!\left(\!\smallmatrix
\ba_{1,r-1}\\ \vspace{1pt}\bh_{s-l,s}^{(1)}\endsmallmatrix\!\!\right) \ot
a_r^{\fh_{s-l,s}^{(2)}}\\
& + \sum_{i=1}^l (-1)^{r+i} F^{(l-i+1)}\!\!
\left(\!\smallmatrix F^{(i)}\!\!\left(\!\smallmatrix\ba\\
\vspace{1pt} \bh_{s-i+1,s}^{(1)} \endsmallmatrix\!\!\right)\\ \vspace{1pt}
\bh_{s-l,s-i}^{(1)}\endsmallmatrix \!\!\!\right) \ot
f(\fh_{s-l,s-i}^{(2)},\fh_{s-i+1,s}^{(2)}).
\endalign
$$
where $F^{(l+1)}\!\!\left(\!\smallmatrix \ba_{1,r-1}\\ \vspace{1pt}
\bh_{s-l,s} \endsmallmatrix\!\!\right) = F^{(l+1)}(\bh_{s-l,s})$ if $r=1$.

\endproclaim

\demo{Proof} We prove the second formula. The proof of the first one is
similar. It is clear that the lemma is valid for $l=1$. Let $l>1$ and
suppose the result is valid for $l-1$. To abbreviate we put
$$
\align
& \xi = u(l-1)+j+s\\
& \fh_{s-l,s}^{j(4)} = \fh_{s-l,j+1}^{(4)}\fh_{j+2,s}^{(3)},\\
& \bh_{s-l,s}^{j(3)} = \bh_{s-l,j-1}^{(3)} \ot h_j^{(3)}h_{j+1}^{(3)} \ot
\bh_{j+2,s}^{(2)},\\
& \bbf_j^{(2)} = f(h_j^{(2)},h_{j+1}^{(2)})^{\ov{\bh_{s-l,j-1}^{(2)}}},\\
& \bbf_{s-l,s-i,s}^{j(4)} = f(\fh_{s-l,j+1}^{(4)}\fh_{j+2,s-i}^{(3)},
\fh_{s-i+1,s}^{(3)}), \\
& \bbf_{s-l,s-i,s}^{(4)j} = f(\fh_{s-l,s-i}^{(4)},\fh_{s-i+1,j+1}^{(4)}
\fh_{j+2,s}^{(3)})\\
& \bbf_{s-l,s-i,s}^{(2)} = f(\fh_{s-l,s-i}^{(2)},\fh_{s-i+1,s}^{(2)}).
\endalign
$$
We have:
$$
\align
& F^{(l+1)}\!\!\left(\!\smallmatrix \ba\\ \bh_{s-l,s} \endsmallmatrix
\!\!\right) = \sum_{j=s-l}^{s-1}\sum_{u=0}^r (-1)^{\xi-1}\ba_{1u}^{
\ov{\bh_{s-l,s}^{(1)}}}\ot\bbf_j^{(2)}\ot F^{(l)}\!\!\left(\!\smallmatrix
\ba_{u+1,r} \\ \vspace{1pt}\bh_{s-l,s}^{j(3)} \endsmallmatrix\!\!\right)\\
& = \sum_{j=s-l}^{s-1} \sum_{u=0}^{r-1} (-1)^{\xi-1}\ba_{1u}^{\ov{
\bh_{s-l,s}^{(1)}}}\ot \bbf_j^{(2)} \ot F^{(l)}\!\!\left(\!
\smallmatrix \ba_{u+1,r-1}\\ \vspace{1pt} \bh_{s-l,s}^{j(3)}
\endsmallmatrix \!\!\right)\ot a_r^{\fh_{s-l,s}^{j(4)}}\\
&+ \sum_{j=s-l}^{s-2}\sum_{u=0}^r\sum_{i=1}^{s-j-1} (-1)^{\xi-1+r-u+i}
\ba_{1u}^{\ov{\bh_{s-l,s}^{(1)}}} \ot \bbf_j^{(2)} \ot F^{(l-i)}
\!\!\left(\!\!\!\smallmatrix F^{(i)}\!\!\left(\!\smallmatrix \ba_{u+1,r}\\
\vspace{1pt}\bh_{s-i+1,s}^{(2)} \endsmallmatrix\!\!\right)\\
\vspace{1pt}\bh_{s-l,s-i}^{j(3)} \endsmallmatrix\!\!\!\right)\ot
\bbf_{s-l,s-i,s}^{j(4)}\\
& + \sum_{j=s-l+1}^{s-1} \sum_{u=0}^r \sum_{i=s-j}^{l-1}
(-1)^{\xi-1+r-u+i} \ba_{1u}^{\ov{\bh_{s-l,s}^{(1)}}}\ot \bbf_j^{(2)} \!
\ot F^{(l-i)} \!\!\left(\!\!\!\smallmatrix F^{(i)}\!\!\left(\!\smallmatrix
\ba_{u+1,r} \\ \vspace{1pt}\bh_{s-i,s}^{j(3)} \endsmallmatrix\!\!\right)\\
\vspace{1pt} \bh_{s-l,s-i-1}^{(3)}\endsmallmatrix\!\!\!\right)\!\ot
\bbf_{s-l,s-i-1,s}^{(4)j}.
\endalign
$$
Permuting the order of the summands, we obtain
$$
\align
& F^{(l+1)}\!\!\left(\!\smallmatrix \ba\\ \bh_{s-l,s} \endsmallmatrix
\!\!\right) = \sum_{j=s-l}^{s-1}\sum_{u=0}^{r-1} (-1)^{\xi-1}
\ba_{1u}^{\ov{\bh_{s-l,s}^{(1)}}}\!\ot \bbf_j^{(2)} \ot F^{(l)} \!\!
\left(\! \smallmatrix \ba_{u+1,r-1}\\ \vspace{1pt} \bh_{s-l,s}^{j(3)}
\endsmallmatrix \!\!\right)\ot a_r^{\fh_{s-l,s}^{j(4)}}\\
&+\sum_{i=1}^{l-1}\sum_{u=0}^r\sum_{j=s-l}^{s-i-1} (-1)^{\xi-1+r-u+i}
\ba_{1u}^{\ov{\bh_{s-l,s}^{(1)}}} \ot \bbf_j^{(2)} \ot F^{(l-i)} \!\!\left(
\!\!\!\smallmatrix F^{(i)}\!\!\left(\!\smallmatrix\ba_{u+1,r}\\\vspace{1pt}
\bh_{s-i+1,s}^{(2)}\endsmallmatrix\!\!\right)\\ \vspace{1pt} \bh_{s-l,s-i}
^{j(3)}\endsmallmatrix\!\!\!\right) \ot \bbf_{s-l,s-i,s}^{j(4)}\\
& + \sum_{i=2}^l \sum_{u=0}^r\sum_{j=s-i+1}^{s-1} (-1)^{\xi+r-u+i}
\ba_{1u}^{\ov{\bh_{s-l,s}^{(1)}}}\ot \bbf_j^{(2)} \ot F^{(l-i+1)}
\!\!\left(\!\!\!\smallmatrix F^{(i-1)} \!\!\left(\!\smallmatrix
\ba_{u+1,r}\\ \vspace{1pt} \bh_{s-i+1,s}^{j(3)}\endsmallmatrix\!\!\right)\\
\vspace{1pt} \bh_{s-l,s-i}^{(3)} \endsmallmatrix\!\!\!\right) \ot
\bbf_{s-l,s-i,s}^{(4)j}\\
& = F^{(l+1)}\!\!\left(\!\smallmatrix \ba_{1,r-1}\\ \vspace{1pt}
\bh_{s-l,s}^{(1)} \endsmallmatrix\!\!\right)\ot a_r^{\fh_{s-l,s}^{(2)}} +
\sum_{i=1}^l (-1)^{r+i} F^{(l-i+1)}\!\!\left(\!\!\!\smallmatrix F^{(i)}
\!\! \left(\!\smallmatrix \ba_{1r}\\ \vspace{1pt} \bh_{s-i+1,s}^{(1)}
\endsmallmatrix \!\!\right)\\ \vspace{1pt} \bh_{s-l,s-i}^{(1)}
\endsmallmatrix\!\!\!\right)\ot \bbf_{s-l,s-i,s}^{(2)},
\endalign
$$
which ends the proof\qed

\enddemo

\subhead Computation of $d^l_{rs}$\endsubhead Let us compute $d_{rs}^{l+1}$
for $l\ge 1$. First we suppose the formula is valid for $d_{rs}^j$ with
$j\le l$ and we see that it is valid for $d_{0s}^{l+1}$. To abbreviate we
write $\zeta_i = is+(l-i+1)(s-1)+1$. Using the inductive hypothesis and the
fact that $\si^0 \circ d^l\bigl(a_0\ot \bh_{0s}\ot 1\# 1\bigr) = 0$, we
obtain:
$$
\align
& d^{l+1}\bigl(1\ot \bh\ot 1_E\bigr) = - \sum_{i=1}^l \si^0
\circ d^{l+1-i} \circ d^i\bigl(1\ot \bh\ot 1_E\bigr)\\
& = \sum_{i=1}^l (-1)^{is+1}\si^0\circ d^{l+1-i}\Bigl(1\ot\bh_{0,s-i}\ot
F^{(i)}(\bh_{s-i+1,s}^{(1)})\ot 1\#\fh_{s-i+1,s}^{(2)}\Bigr)\\
& = \sum_{i=1}^l \si^0\biggl(\! (-1)^{\zeta_i} \! \ot\bh_{0,s-l-1}\! \ot
F^{(l+1-i)} \!\!\left(\!\!\!\smallmatrix F^{(i)}(\bh_{s-i+1,s} ^{(1)})
\vspace{1pt} \bh_{s-l,s-i}^{(1)} \endsmallmatrix\!\!\!\right) \! \ot\!
f(\fh_{s-l,s-i}^{(2)},\fh_{s-i+1,s}^{(2)})\#\fh_{s-l,s}^{(3)}\! \biggr)\\
&= (-1)^{(l+1)s}1\ot \bh_{0,s-l-1}\ot F^{(l+1)}(\bh_{s-l,s}^{(1)})\ot 1\#
\fh_{s-l,s}^{(2)},
\endalign
$$
where the last equality follows from the definition of $\si^0$ and
Lemma~B.3. Now, we suppose the result is valid for $d_{r's}^{l+1}$ with
$r'<r$ and we show that it is valid for $d_{rs}^{l+1}$. To abbreviate we
write $\zeta_i = i(r+s)+(l-i+1)(r+s-1)+1$.
$$
\align
& d^{l+1}\bigl(1\ot \bh\ot \ba\ot 1_E\bigr) = - \sum_{i=0}^l
\si^0 \circ d^{l+1-i} \circ d^i\bigl(1\ot \bh\ot \ba\ot 1_E\bigr)\\
& = (-1)^{r+1} \si_0 \circ d^{l+1}\bigl(1\ot \bh\ot \ba\ot 1\bigr) -
(-1)^{r+s}\si_0\circ d^l\bigl(1\ot \bh_{0,s-1}\ot \ba^{h_s^{(1)}} \ot 1\#
h_s^{(2)}\bigr)\\
& - \sum_{i=2}^l \si^0 \circ d^{l+1-i}\left((-1)^{i(r+s)}\ot \bh_{0,s-i}
\ot F^{(i)} \!\!\left(\!\smallmatrix \ba\\ \vspace{1pt}
\bh_{s-i+1,s}^{(1)}\endsmallmatrix\!\!\right) \ot 1\# \bh_{s-i+1,s}^{(2)}
\right)\\
& = (-1)^{r+1} \si_0 \circ d^{l+1}\bigl(1\ot \bh\ot \ba\ot 1\bigr)\\
& - \sum_{i=1}^l \si^0 \circ d^{l+1-i}\left((-1)^{i(r+s)}\ot \bh_{0,s-i}
\ot F^{(i)}\!\!\left(\!\smallmatrix \ba\\ \vspace{1pt} \bh_{s-i+1,s}^{(1)}
\endsmallmatrix\!\!\right) \ot 1\# \bh_{s-i+1,s}^{(2)} \right)\\
& = \si^0\!\Biggl((-1)^{(l+1)(r+s-1)+r+1}\ot \bh_{0,s-l-1}\ot F^{(l+1)}
\!\!\left(\! \smallmatrix\ba_{1,r-1}\\ \vspace{1pt} \bh_{s-l,s}^{(1)}
\endsmallmatrix\!\!\right)\ot a_r^{\fh_{s-l,s}^{(2)}}\# \fh_{s-l,s}^{(3)}\\
& + \sum_{i=1}^l (-1)^{\zeta_i}\ot \bh_{0,s-l-1} \ot F^{(l+1-i)}
\!\!\left(\!\!\!\smallmatrix F^{(i)} \!\!\left(\!\smallmatrix \ba\\
\vspace{1pt} \bh_{s-i+1,s}^{(1)} \endsmallmatrix\!\!\right)\\ \vspace{1pt}
\bh_{s-l,s-i}^{(1)} \endsmallmatrix \!\!\!\right)\ot f(\bh_{s-l,s-i}^{(2)},
\bh_{s-i+1,s}^{(2)})\# \bh_{s-l,s}^{(3)} \Biggr) \\
& = (-1)^{(l+1)(r+s)}\ot\bh_{0,s-l-1}\ot F^{(l+1)}\!\!\left(\!\smallmatrix
\ba\\ \vspace{1pt}\bh_{s-l,s}^{(1)}\endsmallmatrix\!\!\right)\ot 1\#
\fh_{s-l,s}^{(2)},
\endalign
$$
where the last equality follows from the definition of $\si^0$ and
Lemma~B.3\qed

\remark{Remark B.4} When $H$ is a group algebra $k[G]$ and the
$2$-cocycle $f$ takes its values in the center of $A$, then
$$
d^l_{rs}\bigl(a_0\ot\bg_{0s}\ot\ba_{1r}\ot 1_E\bigr) = (-1)^{l(r+s)} a_0
\ot\bg_{0,s-l}\ot F_0^{(l)}(\bg_{s-l+1,s})*\ba_{1r}\ot 1\# \fg_{s-l+1,s},
$$
where $*$ denotes the shuffle product:
$$
\ba_{1r}*\bb_{1l} = \sum_{0\le i_1\le \dots \le i_r\le l}
(-1)^{i_1+\cdots+i_l} b_1\ot \cdots\ot b_{i_1}\ot a_1\ot b_{i_1+1}\ot
\cdots\ot b_{i_r}\ot a_r\ot b_{i_r+1}\ot \cdots.
$$
\endremark

\enddocument